\documentclass[11pt,twoside]{atmp}

\newtheorem{theorem}{Theorem}[section]
\newtheorem{lemma}[theorem]{Lemma}

\theoremstyle{definition}
\newtheorem{definition}[theorem]{Definition}

\theoremstyle{remark}
\newtheorem{remark}[theorem]{Remark}

\numberwithin{equation}{section}

\input{tcilatex}
\begin{document}

\title[The Conjugate Linearized Ricci Flow...]
{The Conjugate Linearized Ricci Flow on Closed 3--Manifolds}

\author{MAURO CARFORA}



\thanks{Dipartimento di Fisica Nucleare e Teorica, Universita` degli Studi di Pavia\\
and\\
Istituto Nazionale di Fisica Nucleare, Sezione di Pavia\\
via A. Bassi 6, I-27100 Pavia, Italy\\
E-mail address: mauro.carfora@pv.infn.it}




\date{}


\begin{abstract}
We characterize the conjugate linearized Ricci flow and the associated backward heat kernel on closed three--manifolds of bounded geometry. We discuss their  properties, and introduce the notion of Ricci flow conjugated constraint sets which characterizes a way of Ricci flow averaging metric dependent geometrical data. We also provide an integral representation of the Ricci flow metric itself and of its Ricci tensor in terms of the heat kernel of the conjugate linearized Ricci flow. These results, which readily extend to closed $n$--dimensional manifolds, yield for various conservation laws,  monotonicity and asymptotic formulas for the Ricci flow and its linearization.
\end{abstract}
\maketitle

\vfill\eject
\section{INTRODUCTION}

Hamilton's Ricci flow \cite{11} is the weakly--parabolic geometric evolution equation obtained by deforming a Riemannian metric $g_{ab}$, on a smooth $n$--manifold $\Sigma $, in the direction of its Ricci tensor $\mathcal{R}_{ab}$ \cite{1,6,11,12,Huisken}. The geometrical and analytical properties featuring in this natural geometric flow have eventually lead to a remarkable proof, due to G. Perelman \cite{18, 19, 20}, of Thurston's geometrization program for three-manifolds \cite{thurston1,thurston2}. This is a result of vast potential use also in theoretical physics, where the Ricci flow often appears in disguise as a natural real-space renormalization group flow. Non-linear $\sigma $-model theory, describing quantum strings propagating in a background spacetime, affords the standard case study  \cite{Bakas, Bakas2,9, lottCMC, Oliynyk}. Paradigmatical applications occur also in relativistic cosmology \cite{CarfMarz2, CarfKam}, (for a series of recent results see also \cite{4,5,CarFokk} and the references cited therein). An important role both in Ricci flow theory, as well as in its physical applications, is played by its formal linearization around a given Ricci evolving metric $\beta \rightarrow g_{ab}(\beta)$, $0\leq \beta  <T_{0}\leq \infty $. 
By suitably fixing the action of the diffeomorphism group $\mathcal{D}iff(\Sigma )$, this linearized flow takes the form of the parabolic initial value problem
\begin{equation} 
\begin{tabular}{l}
$\frac{\partial  }{\partial  \beta   }{\widetilde h}_{ab}(\beta  )=\Delta _{L}\, {\widetilde h}_{ab}(\beta  )\;,$ \\ 
\\ 
${\widetilde h}_{ab}(\beta  =0)=\,{\widetilde h}_{ab}$\, ,\;\; $0\leq \beta  <T_{0}$\;.%
\end{tabular}
  \label{linDT0}
\end{equation}
where $\Delta _{L}$ denotes the Lichnerowicz--DeRham laplacian \cite{lichnerowicz}, (with respect to $g_{ab}(\beta )$), and the symmetric bilinear form $\beta \mapsto {\widetilde h}_{ab}(\beta )$ can be thought of as representing an infinitesimal deformation, $g^{(t)}_{ab}(\beta )=
g_{ab}(\beta )+ t\;{\widetilde h}_{ab}(\beta )$, $t\in (-\varepsilon ,\varepsilon )$, of the given  flow $\beta \rightarrow g_{ab}(\beta)$.
Stability questions around fixed points of the Ricci flow \cite{cao, guenther, guenther3, lamm, simon, sesum2, Ye}, pinching estimates \cite{01}, and characterization of  linear Harnack inequalities \cite{chowhamilton,chowluni,LeiNi}, are typical issues related to the structure of the linearized Ricci flow (\ref{linDT0}). Related problems, with an impact also in the physical applications of the theory, concerns the control of $\beta \mapsto {\widetilde h}_{ab}(\beta )$ not just around fixed points but along a generic Ricci flow metric $\beta \mapsto g_{ab}(\beta )$. In particular,  if one needs to go beyond a fixed point stability analysis, the characterization of monotonicity formulas for the parabolic equation (\ref{linDT0}) is a key problem in many applications.  Difficulties in dealing with such questions are strictly related to the $\mathcal{D}iff(\Sigma )$--equivariance of the Ricci flow. This remark takes shape in the fact that the flow $\beta \mapsto {\widetilde h}_{ab}(\beta )$, solution of (\ref{linDT0}),  may  describe reparametrization of  $\beta \mapsto g_{ab}(\beta )$ as well as the evolution of non--trivial deformations. The former correspond to the $\mathcal{D}iff(\Sigma )$--solitonic solutions of (\ref{linDT0}). They are provided by ${\widetilde h}_{ab}(\beta )=\mathcal{L}_{v(\beta )}g_{ab}(\beta )$, where $\mathcal{L}_{v(\beta )}$ is the Lie derivative along some suitably chosen $\beta $--dependent vector field $v(\beta )$. 
The latter are instead parametrized by $\beta \mapsto {\widetilde h}_{ab}(\beta )$ with $\nabla ^{a}\,{\widetilde h}_{ab}(\beta )=0$, where the divergence $\nabla \cdot $ is with respect to the $\beta $--varying Ricci flow metric $g_{ab}(\beta )$. As is well known,
the  subspace generated by the Lie derivative of the metric along  vector fields,  and the subspace of divergence--free ${\widetilde h}_{ab}$\,'s, provide an $L^{2}(\Sigma ,g)$--orthogonal splitting of the whole space of symmetric bilinear forms. It is a matter of fact, naturally related to the geometry of the Ricci flow, that  (\ref{linDT0}) does not preserve such a splitting unless the Ricci flow $\beta \mapsto g_{ab}(\beta )$ is restricted to particular class of geometries \cite{avez, buzzanca, guenther, guenther3, Ye}.
What happens is that (\ref{linDT0}) may evolve a divergence--free ${\widetilde h}_{ab}(\beta=0)$ into a flow 
$\beta \mapsto {\widetilde h}_{ab}(\beta)$ possessing also a Lie--derivative part. For instance, if one considers, for the volume--normalized Ricci flow, the evolution of the coupled $\beta \mapsto (g_{ab}(\beta ),\,{h}_{ab}(\beta ))$ with $Ric(g)|_{\beta =0}>0$, $\Sigma \simeq \mathbb{S}^{3}$, then by Hamilton's rounding theorem \cite{11}, $g_{ab}(\beta )$ converges, as $\beta \nearrow \infty $, to the standard metric $\bar{g}$ on the $3$--sphere $\mathbb{S}^{3}$, with $Vol\,[\mathbb{S}^{3},\,\bar{g}]=Vol\,[\mathbb{S}^{3},\,{g(\beta =0)}]$. Since $(\mathbb{S}^{3},\,\bar{g})$ is isolated (\emph{i.e.}, it does not admit any non--trivial Riemannian deformation), it follows that any divergence--free ${\widetilde h}_{ab}(\beta=0))$
 must necessarily evolve under (the normalized version of) (\ref{linDT0}) into a Lie derivative term $\mathcal{L}_{X}\,\bar{g}_{ab}$, for some $\beta $--dependent vector field $X$. This dynamical generation of $\mathcal{D}iff(\Sigma )$--reparametrization out of  non--trivial deformations is at the root of the difficulties in the general  analysis of  (\ref{linDT0}). A possible way out is to adopt a strategy akin to the one used by Perelman \cite{18} in handling Ricci flow $\mathcal{D}iff(\Sigma )$--solitons. These are put under control by means of a (backward) diffusion process which is conjugated to the Ricci diffusion of the Riemannian measure. A related and very subtle use of the backward--forward conjugation, in connection with the K$\ddot a$hler-Ricci flow, has been recently pointed out also by Lei Ni \cite{Lei1}. By extending these points of view to the $\mathcal{L}_{v(\beta )}g_{ab}(\beta )$ solitonic solutions of (\ref{linDT0}),  we introduce in this paper the  backward conjugated flow associated with (\ref{linDT0}), generated by the operator
\begin{equation}
\bigcirc^{*} _{L}\doteq -\frac{\partial }{\partial \beta }-\triangle _{L} +\mathcal{R}\;,
\label{conjflow0}
\end{equation}
where $\mathcal{R}$ denotes the scalar curvature of  $(\Sigma ,g_{ab}(\beta ))$.
The  flow described by $\bigcirc^{*} _{L}$ enjoys many significant properties:

\noindent \emph{(i)} The space of divergence--free bilinear forms is an invariant subspace of the flow.

\noindent \emph{(ii)} If $\beta \mapsto {\widetilde h}_{ab}(\beta )$ is the  a solution of the linearized Ricci flow (\ref{linDT0}), and $\eta \mapsto H^{ab}_{(T)}(\eta )$, $\eta\doteq \beta ^{*}-\beta $, for some $\beta ^{*}\in (0,T_0)$, is a divergence--free solution of the conjugate flow, $\bigcirc^{*} _{L}\,H^{ab}_{(T)}(\eta )=0$, then 
\begin{equation}
\int_{\Sigma }{\widetilde h}_{ab}^{(T)}(\eta )\,\,H^{ab}_{(T)}(\eta )\,d\mu _{g(\eta )}\;,
\end{equation}
where ${\widetilde h}_{ab}^{(T)}(\eta )$ is the divergence--free part of ${\widetilde h}_{ab}(\eta )$, is a conserved quantity along the (backward) Ricci flow.
This result provides a useful control on the dynamics of $\beta \mapsto {\widetilde h}_{ab}^{(T)}(\beta )$.

\noindent \emph{(iii)} If $\beta \mapsto \mathcal{R}_{ab}(\beta )$ is the Ricci tensor associated with the Ricci flow metric $\beta \mapsto g_{ab}(\beta )$, and $\bigcirc^{*} _{L}\,H^{ab}(\eta )=0$, then
\begin{equation}
\int_{\Sigma }R_{ab}(\eta )H^{ab}(\eta )d\mu _{g(\eta )}\;,
\end{equation}
and
\begin{equation}
\int_{\Sigma }\left(g_{ab}(\eta)-2\eta\,R_{ab}(\eta )\right)H^{ab}(\eta )d\mu _{g(\eta )}\;,
\end{equation}
are also conserved along the (backward) Ricci flow. Thus, quite surprisingly, the conjugate linearized Ricci flow has strong averaging properties on the full Ricci flow itself. These averaging properties become manifest when we identify the flow $\eta \mapsto H^{ab}(\eta )$,  with the (backward) heat kernel of $\bigcirc^{*} _{L}$. In such a setting we prove the main results of this paper, \emph{viz.},

\begin{proposition}
Let $\eta \mapsto g_{ab}(\eta )$ be a backward Ricci flow with bounded geometry on $\Sigma_{\eta }\times [0,\beta ^{*}]$ and let ${K}^{ab}_{i'k'}(y,x;\eta )$ be the (backward) heat kernel of the corresponding conjugate linearized Ricci operator 
$\bigcirc ^{*}_{L}\, {K}^{ab}_{i'k'}(y,x;\eta )=0$, for  $\eta \in (0,\beta ^{*}]$, \, with ${K}^{ab}_{i'k'}(y,x;\eta\searrow 0^{+} )={\delta }^{ab}_{i'k'}(y,x)$, (the bi--tensorial Dirac measure). Then 
\begin{equation}
\mathcal{R}_{i'k'}(y,\eta=0 )=\int_{\Sigma }{K}^{ab}_{i'k'}(y,x;\eta )\,\mathcal{R}_{ab}(x,\eta )\,d\mu _{g(x,\eta )}
\;,
\end{equation}
for all $0\leq \eta \leq \beta ^{*}$.
Moreover, as $\eta\searrow 0^{+} $, we have the uniform asymptotic expansion
\begin{eqnarray}
&&\;\;\;\;\;\mathcal{R}_{i'k'}(y,\eta=0 )=\\
\nonumber\\
&&\frac{1}{\left(4\pi \,\eta \right)^{\frac{3}{2}}}\,\int_{\Sigma }\exp\left(-\frac{d^{2}_{0}(y,x)}{4\eta } \right)\,{\tau }^{ab}_{i'k'}(y,x;\eta )\,\mathcal{R}_{ab}(x,\eta )\,d\mu _{g(x,\eta )}\nonumber\\
\nonumber\\
&&+\sum_{h=1}^{N}\frac{\eta ^{h}}{\left(4\pi \,\eta \right)^{\frac{3}{2}}}\,\int_{\Sigma }\exp\left(-\frac{d^{2}_{0}(y,x)}{4\eta } \right)\,{\Phi [h] }^{ab}_{i'k'}(y,x;\eta )\,\mathcal{R}_{ab}(x,\eta )\,d\mu _{g(x,\eta )}\nonumber\\
\nonumber\\
&&+O\left(\eta ^{N-\frac{1}{2}} \right)\nonumber\;,
\end{eqnarray}
\noindent where ${\tau }^{ab}_{i'k'}(y,x;\eta )$ $\in T\Sigma_{\eta }\boxtimes T^{*}\Sigma_{\eta } $  is the parallel transport operator
associated with  $(\Sigma ,g(\eta ))$,   $d_{0}(y,x)$ is the distance function in $(\Sigma ,g(\eta=0 ))$, and ${\Phi [h] }^{ab}_{i'k'}(y,x;\eta )$ are the smooth section  $\in C^{\infty }(\Sigma\times \Sigma ' ,\otimes ^{2}T\Sigma\boxtimes \otimes ^{2}T^{*}\Sigma)$, (depending on the geometry of $(\Sigma ,g(\eta ))$), characterizing the asymptotics of the heat kernel ${K}^{ab}_{i'k'}(y,x;\eta )$.
\label{princ1}
\end{proposition}
\noindent Under the same hypotheses of proposition \ref{princ1}, we also have the following integral representation of the Ricci flow on $\Sigma _{\beta }\times (0,\beta ^{*}]$.

\begin{proposition}
Let $\beta \mapsto g_{ab}(\beta )$ be a Ricci flow with bounded geometry on $\Sigma_{\beta }\times [0,\beta ^{*}]$, and let ${K}^{ab}_{i'k'}(y,x;\eta )$ be the (backward) heat kernel of the corresponding conjugate linearized Ricci operator 
$\bigcirc ^{*}_{L}$, for  $\eta=\beta ^{*}-\beta $. Then, along the backward flow $\eta \mapsto g_{ab}(\eta )$,
\begin{equation}
 g_{i'k'}\,(y,\eta=0)=\int_{\Sigma }{K}^{ab}_{i'k'}(y,x;\eta)\,\left[{g}_{ab}(x,\eta)-2\eta \,\,\mathcal{R}_{ab}(x,\eta)\right]\,d\mu _{g(x,\eta )}\;,
\label{grepres0}
\end{equation}
\noindent for all $0\leq \eta \leq \beta ^{*}$, and 
\begin{eqnarray}
&&\;\;\;\;\;\;\;\;{g}_{i'k'}(y,\eta=0)=\\
\nonumber\\
&&\frac{1}{\left(4\pi \,\eta \right)^{\frac{3}{2}}}\,\int_{\Sigma }e^{-\frac{d^{2}_{0}(y,x)}{4\eta }}\,{\tau }^{ab}_{i'k'}(y,x;\eta )\,\left[{g}_{ab}(x,\eta)-2\eta\,\mathcal{R}_{ab}(x,\eta)\right]\,d\mu _{g(x,\eta)}     \nonumber\\
\nonumber\\
&&+\sum_{h=1}^{N}\frac{\eta ^{h}}{\left(4\pi \,\eta \right)^{\frac{3}{2}}}\,\int_{\Sigma }e^{-\frac{d^{2}_{0}(y,x)}{4\eta }}\,{\Phi [h] }^{ab}_{i'k'}(y,x;\eta )\left[{g}_{ab}(x,\eta)-2\eta\,\mathcal{R}_{ab}(x,\eta)\right]\,d\mu _{g(x,\eta)} \nonumber\\
\nonumber\\
&&+O\left(\eta ^{N-\frac{1}{2}} \right)\nonumber\;.
\end{eqnarray}
holds uniformly, as $\eta \searrow 0^{+}$.
\end{proposition}
\noindent In particular, the above result proves the following 
\begin{theorem}
The heat kernel flow
\begin{equation}
\eta\longmapsto {K}^{ab}_{i'k'}(y,x;\eta)\;
\end{equation}
is conjugated and thus fully equivalent to the Ricci flow $\beta\longmapsto g_{ab}(\beta)$.
\end{theorem}
\noindent This is a quite non--trivial consequence of the conjugacy relation and opens the possibility of a weak formulation of the Ricci flow by exploiting the linear evolution of $\eta\longmapsto {K}^{ab}_{i'k'}(y,x;\eta)$.

\noindent The properties of the conjugate heat flow \cite{ecker},\cite{LeiNi},\cite{18} and those of the conjugate linearized Ricci flow established in this paper  suggest to shift emphasis from the flows themselves to their dependence from the corresponding initial data. Thus, along a Ricci flow of bounded geometry $\beta \mapsto (\Sigma ,{g}(\beta ))$, $\beta\in [0,\beta ^{*}]$, we consider the associated heat flow $(\beta,\varrho _{0}) \mapsto \varrho (\beta )$, $(\frac{\partial}{\partial\beta}-\Delta )\,\varrho=0$, and linearized Ricci flow $(\beta,h_{ab}) \mapsto {\widetilde{h}}_{ab}(\beta )$, as functionals of the respective initial data  $\varrho (\beta=0)\doteq \varrho_{0}$, and ${\widetilde h}_{ab}(\beta =0)\doteq h_{ab}$. Similarly, we can consider, along the backward Ricci flow 
$\eta \mapsto (\Sigma ,{g}(\eta ))$, $\eta\in [0,\beta ^{*}]$, $\eta\doteq \beta^{*}-\beta$, the conjugate flows $(\eta,\varpi _{*}) \mapsto \varpi (\eta )$, $(\frac{\partial}{\partial\eta}-\Delta+\mathcal{R} )\,\varpi=0$,  and $(\eta,H_{*}^{ab}) \mapsto {H}^{ab}(\eta )$, as functionals of the respective initial data  $\varpi (\eta=0)\doteq \varpi_{*}$, and $H^{ab}(\eta =0)\doteq H_{*}^{ab}$. In general, if the initial data $(\varrho_{0},h_{ab})$ satisfy a geometrical condition in the form of a constraint $\mathcal{C}(g_{ik}(\beta=0),\varrho_{0},h_{ab})=0$, then this constraint will not be preserved along the evolution of the given data. However, if we are able to find, along the given Ricci flow, initial data 
$(\varpi_{*},H^{ab}_{*})$ for the conjugated flows such that $\mathcal{C}(g_{ik}(\eta=0),\varpi_{*},H^{ab}_{*})=0$, then
the conjugate flows interpolate between $(\varrho_{0},h_{ab})$ and $(\varpi_{*},H^{ab}_{*})$ by averaging the data with the \emph{kernels} $(\varpi(\eta),H^{ab}(\eta))$. We say, in such a case, that the constraints $\mathcal{C}(g_{ik}(\beta=0),\varrho_{0},h_{ab})=0$ and $\mathcal{C}(g_{ik}(\eta=0),\varpi_{*},H^{ab}_{*})=0$ are Ricci flow conjugated. This is basically a way for averaging geometrical constraints along the Ricci flow, and may find applications in various geometrical and physical setting. The stability of Type--$II$ singularities (see Section 7), and
the problem of Ricci flow deforming the initial data set for Einstein equations (see Section 4) may provide important examples.

Coming to the structure of the paper, we have tried to keep the presentation as self--contained as possible. We start by recalling some  well--known facts about the Ricci flow and its linearization in Section 2. The conjugate linearized Ricci flow is introduced in Section 3, where we also establish its properties.  In Section 4 we discuss  the  heat kernel associated with the conjugate linearized Ricci flow. In an appendix, kindly provided by Stefano Romano, we carry out the explicit construction of the the heat kernel of a generalized Laplacian when the operator in question smoothly depend on a one--parameter family of metrics. The results discussed in Section 5 are elementary consequences of the properties of the heat kernel of the conjugate linearized Ricci flow. In Section 6 we formalize the notion of Ricci flow conjugated constraint sets, and briefly discuss a few natural examples.

\subsection*{Aknowledgements}
The author would like to thank Giovanni Bellettini, Carlo Mantegazza and Matteo Novaga for useful conversations in the preliminary stage of preparation of this paper.   Gerhardt Huisken deserves a special mention for his kind suggestions and comments on the possibility of applying the backward linearized Ricci flow in studying stability for Type--$II$ singularities. A special debt of gratitude goes to Stefano Romano for his invaluable help in the heat kernel analysis, to Lei Ni for his comments and many suggestions, and to my dear friend Thomas Buchert for many stimulating ideas that had an important influence on the problems discussed here.\\
\noindent The research of this paper has been partially supported by the 2006-PRIN Grant \emph{Spazi di moduli, strutture integrabili, e equazioni alle derivate parziali di tipo evolutivo}.

\section{\protect\medskip Remarks on the  Ricci flow and its linearization}
We start by collecting a number of  technical results on Ricci flow theory that we shall need in the sequel. Excellent sources of information are provided by \cite{1}, \cite{6}, \cite{Glickenstein}, \cite{chowluni}, and \cite{topping}. For simplicity $\Sigma$ will always denote a $C^{\infty }$ compact three-dimensional manifold without boundary, and $C^{\infty }(\Sigma, \mathbb{R})$ and 
${C}^{\infty }(\Sigma,\otimes ^{p}\, T^{*}\Sigma\otimes ^{q}T\Sigma )$ are the space of smooth functions and of smooth $(p,q)$--tensor fields over $\Sigma $, respectively.  
We shall denote by $\mathcal{D}iff(\Sigma )$ the group of smooth diffeomorphisms of $\Sigma $, and by $\mathcal{M}et(\Sigma )$ the space of all smooth Riemannian metrics over $\Sigma$. The tangent space\,, $\mathcal{T}_{(\Sigma ,g)}\mathcal{M}et(\Sigma )$, to $\mathcal{M}et(\Sigma )$ at $(\Sigma,g )$ can be naturally identified with the space of symmetric bilinear forms ${C}^{\infty }(\Sigma,\otimes ^{2}_{S}\, T^{*}\Sigma)$ over $\Sigma $, endowed with the pre--Hilbertian  $L^{2}$ inner product  $(U,V)_{L^{2}(\Sigma )}\doteq \int_{\Sigma }g^{il}\,g^{km}\,U_{ik}\,V_{lm}d\mu _{g}$ for  $U,\,\,V\,\in\,{C}^{\infty }(\Sigma,\otimes ^{2}_{S}\, T^{*}\Sigma)$. 
Let ${L}^{2}(\Sigma,\otimes ^{2}\, T^{*}\Sigma)$  be the corresponding $L^{2}$ completions of ${C}^{\infty }(\Sigma,\otimes ^{2}_{S}\, T^{*}\Sigma)$. A geometric property of $\mathcal{M}et(\Sigma )$ that we shall often exploit is that the tangent space $\mathcal{T}_{(\Sigma ,g)}\mathcal{O}_{g}$ to the $\mathcal{D}iff(\Sigma )$--orbit of a given metric $g\in \mathcal{M}et(\Sigma )$ is the image of the injective operator
\begin{eqnarray}
 \delta_{g} ^{*} \,: &C^{\infty }(\Sigma ,T^{*}\,\Sigma )& \longrightarrow C^{\infty }(\Sigma ,\otimes ^{2}T^{*}\,\Sigma )\label{menoLie}\\
& w_{a}\,dx^{a}& \longmapsto \delta_{g}^{*}\,(w_{a}\,dx^{a})\doteq \frac{1}{2}\,\mathcal{L}_{w^{\#}}\,g\notag\;,
\end{eqnarray}
where we have set $(w^{\#})^{i}\doteq g^{ik}w_{k}$, and denoted by $\mathcal{L}_{w^{\#}}$ the corresponding Lie derivative. Standard elliptic theory then implies that the $L^{2}$--orthogonal subspace to $Im\;\delta_{g} ^{*}$ in $\mathcal{T}_{(\Sigma ,g)}\mathcal{M}et(\Sigma )$ is spanned by the ($\infty $--dim) kernel of the $L^{2}$ adjoint $\delta_{g}$ of  $\delta_{g} ^{*}$,
\begin{eqnarray}
\delta_{g} \,: &C^{\infty }(\Sigma ,\otimes ^{2}T^{*}\,\Sigma )& \longrightarrow C^{\infty }(\Sigma ,T^{*}\,\Sigma )\label{menodiv}\\
 & h_{ab}\,dx^{a}\otimes dx^{b} & \longmapsto \delta_{g} \,(h_{ab}\,dx^{a}\otimes dx^{b})\doteq -\,g^{ij}\,\nabla _{i}h_{jk}\,dx^{k}\;.\nonumber
\end{eqnarray}
It follows that with respect to the inner product $(\circ ,\circ )_{L^{2}(\Sigma )}$, the tangent space $\mathcal{T}_{(\Sigma ,g)}\mathcal{M}et(\Sigma )$ splits as \cite{ebin}
\begin{equation}
\mathcal{T}_{(\Sigma ,g)}\mathcal{M}et(\Sigma )\cong  Ker\;\delta_{g} \,\oplus Im\;\delta_{g} ^{*}\;.
\label{L0split}
\end{equation}
Unless $\mathcal{R}ic(g)\equiv C\,g+\mathcal{L}_{w^{\#}}\,g$, for some constant $C$, the Ricci tensor $\mathcal{R}ic(g)$ of a metric $g\in \mathcal{M}et(\Sigma )$ can be thought of as a non--trivial $\mathcal{D}iff(\Sigma )$--equivariant section of the tangent bundle $\mathcal{T}\,\mathcal{M}et(\Sigma )$, \emph{i.e.}, $\{\mathcal{R}ic(g)\}\cap Ker\;\delta_{g}\not=\emptyset $. Thus, according to (\ref{L0split}), the Ricci flow associated with a Riemannian three-manifold $(\Sigma,g)$ can be thought of as the dynamical system on $\mathcal{M}et(\Sigma )$ generated by the weakly-parabolic diffusion--reaction PDE \cite{11}
\begin{equation} 
\begin{tabular}{l}
$\frac{\partial }{\partial \beta }g_{ab}(\beta )=-2\mathcal{R}_{ab}(\beta ),$ \\ 
\\ 
$g_{ab}(\beta =0)=g_{ab}$\, ,\;\; $0\leq \beta <T_{0}$\;,%
\end{tabular}
   \label{mflow}
\end{equation}
where $\mathcal{R}_{ab}(\beta )$ is the Ricci tensor of the metric $g_{ik}(\beta )$. The flow $(\Sigma ,g) \mapsto (\Sigma ,g(\beta ))$, defined by (\ref{mflow}), always exists in a maximal interval $0\leq \beta \leq T_{0}$,  for some $T_{0}\leq \infty $. If such a $T_{0}$ is finite then 
$\lim_{\beta \nearrow T_{0}}\, [\sup_{x\in \Sigma }\,|Rm(x,\beta )|]=\infty $, \cite{11,12} where $Rm(\beta )$ is the Riemann tensor of $(\Sigma ,g(\beta ))$. Note that, by exploiting  a  result by N. Sesum and M. Simon\cite{sesum, simon2}, (see also the comments in \cite{knopfsimon})
the curvature singularity regime for the $3$--d Ricci flow is equivalent to $\limsup_{\beta \nearrow T_{0}}\, [\max_{x\in \Sigma }\,|Ric(x,\beta )|]=\infty$, (quite surprisingly, this result of Sesum and Simon holds on any compact $n$--dimensional manifold).
The structure of singularities of the Ricci flow, as well as that of generalized fixed points attained if $T_0=\infty$, is associated with self--similar solutions generated by the action of $\mathcal{D}iff(\Sigma )\times \mathbb{R}_{+}$, where $\mathbb{R}_{+}$ acts by scalings. These solutions are described by the Ricci solitons $-2{\mathcal{R}}_{ab}(\beta )=\mathcal{L}_{\vec{v}(\beta )}\,{g}_{ab}+\, \varepsilon \,{g}_{ab}$, where  $\mathcal{L}_{\vec{v}(\beta )}$ denotes the Lie derivative along the $\beta $--dependent (complete) vector field ${\vec{v}(\beta )}$ generating $\beta \mapsto \varphi (\beta )$\,$\in\,\mathcal{D}iff(\Sigma )\times \mathbb{R}_{+}$, and where, up to rescaling, we may assume that $\varepsilon =\,-1,\; 0,\; 1$, (respectively yielding for the shrinking, steady, and expanding solitons). This non--trivial action of the diffeomorphisms group 
$\mathcal{D}iff(\Sigma )$ on the evolution of $(\Sigma ,g(\beta ))$ can be better seen if we describe the kinematics of the flow (\ref{mflow}) in the \emph{parabolic spacetime} $M^{4}_{Par}\doteq \Sigma \times I$,\; $I\doteq [0,T_{0})\subset \mathbb{R}$. We assume that the diffeomorphism 
\begin{equation}
F_{\beta }:I\times \Sigma \longrightarrow M^{4}_{Par};\;\;\; (\beta ,x)\mapsto i_{\beta }(x)\;,
\label{Fdiff}
\end{equation} 
of $I\times \Sigma$ onto $M^{4}_{Par}$, is the identity map, and that $M^{4}_{Par}$ carries the product metric ${}^{(4)}g_{par}$, so that in the coordinates induced by $F_{\beta }$  we can write 
\begin{equation}
(F_{\beta }^{*}\;{}^{(4)}g_{par})=g_{ab}(\beta )dx^{a}\otimes dx^{b}+d\beta \otimes d\beta \;.
\label{Fmetr}
\end{equation}
In such a framework, $\frac{\partial }{\partial \beta }:\Sigma \rightarrow TM^{4}_{Par}$, can be interpreted as a vector field, transversal (actually, ${}^{(4)}g_{par}$--normal) to the leaves $\{\Sigma _{\beta }\}$, describing the Ricci flow evolution as seen by observers at rest on $\Sigma _{\beta }$. The  evolution of the metric $g(\beta )$ can be equivalently described by observers in \emph{motion} on $\Sigma _{\beta }$. To this end,  consider a curve of diffeomorphisms $I\ni \beta\mapsto \varphi (\beta )\in \mathcal{D}iff(\Sigma ) $, (with the initial condition $\varphi ^{i}(x^{a},\beta =0)=id_{\Sigma }$), and define the vector field $X_{\varphi }:\Sigma _{\beta }\rightarrow T\Sigma _{\beta }$, $X_{\varphi }=\frac{\partial }{\partial \beta }\varphi (\beta )$, generating $\beta\mapsto \varphi (\beta )$. Such a $\beta $--dependent  $X_{\varphi }$ provides the velocity field of these non--static observers. Thus,
\begin{equation}
\frac{d}{d\beta }\,F_{\beta,\varphi  }=\frac{\partial }{\partial \beta }+X_{\varphi }\,:\Sigma _{\beta }\longrightarrow T\,M^{4}_{Par},
\label{Rtime}
\end{equation}
is the space--time vector field covering the diffeomorphism  $F_{\beta,\varphi  }$ of $I\times \Sigma$ onto $(M^{4}_{Par}, {}^{(4)}g_{par})$, defining space--time coordinates $(\beta ,y^{i}=\varphi^{i} (\beta ,x))$ which describe the curve of embeddings $(\beta , x)\hookrightarrow (\beta ,\varphi (\beta ,x))$ of $\Sigma _{\beta }$ in $M^{4}_{Par}$. In terms of the coordinates $(\beta ,y^{i})$ we can write
\begin{equation}
(F_{\beta,\varphi   }^{*}\;{}^{(4)}g_{par})=\breve {g}_{ab}(\beta  )(dy^{a}+X^{a}_{\varphi }d\beta )\otimes (dy^{b}+X^{b}_{\varphi }d\beta )+d\beta  \otimes d\beta  \;,
\end{equation}
where the metric $\breve {g}_{ab}(y^{i},\beta)$ is $\beta$--propagated according to 
the Hamilton--DeTurck flow
\begin{equation} 
\begin{tabular}{l}
$\frac{\partial }{\partial \beta  }\breve {g}_{ab}(\beta  )=-2\breve {\mathcal{R}}_{ab}(\beta  )%
-\mathcal{L}_{X_{\varphi }}\breve {g}_{ab}(\beta ),$ \\ 
\\ 
$\breve {g}_{ab}(\beta  =0)=g_{ab}$\, ,\;\; $0\leq \beta  <T_{0}$\;.%
\end{tabular}
  \label{mflowDT}
\end{equation}
The connection between (\ref{mflowDT}) and (\ref{mflow}) is most easily established if we proceed as in the mechanics of continuous media, when shifting from the \emph{body} (Lagrangian) to the  \emph{space} (Eulerian) point of view. To this end, let us introduce the \emph{substantial} derivative $\frac{D}{D\beta }\doteq \frac{\partial }{\partial \beta }+\mathcal{L}_{X_{\varphi }}$ associated with the convective action defined by  $X_{\varphi }$. Since 
$\frac{D }{D \beta  }\,\breve {g}_{ab}(\beta  )=(\varphi ^{*})^{-1}\;\frac{\partial }{\partial \beta }\left[\varphi ^{*}\,\breve {g}\right]_{ab}$ and $\varphi ^{*}\,\mathcal{R}ic(\breve{g})=\mathcal{R}ic(\varphi ^{*}\,\breve {g})$, $\mathcal{R}(\breve {g})=\mathcal{R}(\varphi ^{*}\,\breve {g})$, it follows from (\ref{mflowDT}) that the pull--back $\beta \mapsto (\varphi^{*}\,\breve {g})$ of the flow $\beta \mapsto \breve {g}_{ik}dy^{i}\otimes dy^{k}$, under the action of the $\beta  $-- dependent diffeomorphism  $x^{a}\mapsto y^{i}=\varphi ^{i}(x^{a},\beta)$, solves (\ref{mflow}). The non--trivial action of the diffeomorphism group described above is the rationale underlying
DeTurck's technique for fixing a gauge $F_{\eta,\varphi   }$ making the evolution $\beta \mapsto (T\Sigma, g_{ab}(\beta ))$ of the metric in the tangent bundle  manifestly parabolic \cite{DeTurck}. In this connection, one easy but useful information is that, along the evolution $\beta \mapsto (T\Sigma, g_{ab}(\beta ))$, we can also consider a $\beta$--dependent isomorphism $\iota(\beta )$ between a fixed vector bundle $V$ over $\Sigma $ and the tangent bundle $(T\Sigma , g(\beta))$, in such a way that $(V,\iota^{*}(\beta )\,g(\beta ))$ is isometric to $(T\Sigma , g(\beta=0))$.  This is the Uhlenbeck trick \cite{11a, Ivey}. We briefly describe it to set notation for later use, (what follows holds for any dimension $n$). Consider a bundle isometry
between a fixed vector bundle $V$ over $\Sigma $ and the tangent bundle $T\Sigma $,\,\,
${\iota_{(0)} }\;:\,\left(V, {\iota_{(0)}}^*\,g\right)\longrightarrow \left(T\Sigma ,g \right)$,
where $g$ is a given metric on $T\Sigma $. Locally, in any open set $U\subset \Sigma $, given a basis of sections  $\left\{{e}_{(\mu  )}\right\}_{\mu =1,2,3}$ of $V|_{U}$, and a basis of sections $\{\theta ^{(\nu )}\}_{\nu =1,2,3}$ of the dual bundle $V^{*}|_{U}$, we can write ${\iota_{(0)} }^*\,g|_{U}= {\iota_{(0)}} ^{h}_{\mu }\,\,{\iota_{(0)}} ^{k}_{\nu }\,\,g_{hk}\,\theta ^{(\mu) }\,\theta ^{(\nu) }$, where the components ${\iota_{(0)}} ^{h}_{\mu }$ of the the bundle isomorphism $\iota_{(0)}$ are defined by $\iota _{(0)}(e_{(\mu)})={\iota_{(0)}} ^{h}_{\mu }\,\partial _{h}$.  Let us evolve the isometry $\iota_{(0)} $, along with the Ricci flow $\beta \mapsto g_{ab}(\beta )$,\ $0\leq \beta <T_{0}$, according to $\iota (\beta ):\,(V,{\iota_{(0)}}^*\,g)\rightarrow (T\Sigma , g(\beta ))$, where $\beta \mapsto \iota(\beta) $ is the solution of
\begin{equation} 
\begin{tabular}{l}
$\frac{\partial  }{\partial  \beta  }{\iota ^{k}_{\mu }}(\beta  )=\,{\iota ^{h}_{\mu }}(\beta)\,\, {\mathcal{R}}_{h}^{k}(\beta  ),$ \\ 
\\ 
$\;\;{\iota ^{h}_{\mu }}(\beta  =0)={\iota_{(0)}} ^{h}_{\mu }$\, ,\;\; $0\leq \beta  <T_{0}$\;.%
\end{tabular}
   \label{flowHulen}
\end{equation}
It is easily checked that along such an evolution we have 
$\left( {\iota (\beta ) }^*\,g(\beta )\right)_{\mu\nu}=\left( {\iota_{(0)} }^*\,g\right)_{\mu\nu}$, $0\leq \beta  <T_{0}$, as required. One can also pull back to $(V,{\iota (\beta ) }^*\,g(\beta ))$ the Levi--Civita connection $\nabla(\beta )$ on $(T\Sigma ,g(\beta ))$ according to 
\begin{eqnarray}
 D(\beta ):&& C^{\infty }(\Sigma ,T\Sigma )\times C^{\infty }(\Sigma ,V)\rightarrow  C^{\infty }(\Sigma ,V)\label{iconn}\\
&& (X,\xi ) \longmapsto  D(\beta )_{X}\,\xi:=\iota(\beta )^{*}\,{\nabla(\beta )}_{X}\,\xi \nonumber\;.
\end{eqnarray}
From the defining relation relation $\iota\,[D(\beta )_{h}\,\xi]=\nabla(\beta)_{h}(\iota(\xi))$ it follows that $\nabla_{k}\nabla_{h}(\iota(\xi))$ $=$ $\nabla_{k}\left\{\iota\,[D(\beta )_{h}\,\xi]\right\}$ $=$ $\left\{\iota\,[D(\beta)_{k}D(\beta )_{h}\,\xi]\right\}$. Thus, $\triangle (\beta)(\iota(\xi))$ $:=$ $g^{kh}(\beta)\nabla_{k}\nabla_{h}(\iota(\xi))$ $=$ $g^{kh}(\beta)\left\{\iota\,[D(\beta)_{k}D(\beta )_{h}\,\xi]\right\}$
$=$ $\iota\,[ \Delta _{D}(\beta)\,\xi]$, where the (rough) Laplacian $\Delta _{D}$ on  $(V,{\iota (\beta ) }^*\,g(\beta ))$ is defined by 
\begin{equation}
 \Delta _{D}(\beta)\,\xi:=\,g^{kh}(\beta)\,D(\beta)_{k}D(\beta )_{h}\,\xi\;.
\end{equation}
These remarks imply  a well--known result (see \emph{e.g.},\cite{6},\cite{topping})  that can be phrased in the following form, more adapted to our purposes,  
\begin{lemma}
\label{uhlenlemma}
If a bilinear form $v_{ik}\in {C}^{\infty }(\Sigma,\otimes ^{2}_{S}\, T^{*}\Sigma)$  evolves, along a given Ricci flow $\beta \mapsto g_{ab}(\beta )$,\ $0\leq \beta <T_{0}$, according to the solution  $\beta \mapsto v_{ik}(\beta) $ of the parabolic initial value problem 
\begin{equation} 
\begin{tabular}{l}
$\frac{\partial  }{\partial  \beta  }{v _{ik}}(\beta  )=\,\Delta v_{ik} (\beta),$ \\ 
\\ 
$\;\;{v_{ik}}(\beta  =0)={v_{ik}}$\, ,\;\; $0\leq \beta  <T_{0}$\;,%
\end{tabular}
 .  \label{flowHu}
\end{equation}
then its pull--back $\iota^{*}\,v$, under the map $\iota (\beta ):\,(V,{\iota_{(0)}}^*\,g)\rightarrow (T\Sigma , g(\beta ))$, evolves according to
\begin{equation} 
\begin{tabular}{l}
$\frac{\partial  }{\partial  \beta  }\left({\iota^{i}_{\mu }(\beta )\,\iota^{k}_{\nu }(\beta )\,v _{ik}}(\beta  )\right)=\,\Delta_{D} \left({\iota^{i}_{\mu }(\beta )\,\iota^{k}_{\nu }(\beta )\,v _{ik}}(\beta  )\right)$ \\ 
\\ 
$+{\iota^{h}_{\mu }(\beta )\,\iota^{k}_{\nu }(\beta )\,{\mathcal{R}}_{h}^{i}(\beta  )\,v _{ik}}(\beta  )+
{\iota^{i}_{\mu }(\beta )\,\iota^{h}_{\nu }(\beta )\,{\mathcal{R}}_{h}^{k}(\beta  )\,v _{ik}}(\beta  )$\;,\\
\\
$\left({\iota^{i}_{\mu }(\beta )\,\iota^{k}_{\nu }(\beta )\,v _{ik}}(\beta  )\right)_{(\beta  =0)}=
{{\iota _{(0)}}^{i}_{\mu }\,{\iota_{(0)}}^{k}_{\nu }\,v _{ik}}$\, ,\;\; $0\leq \beta  <T_{0}$\;.%
\end{tabular}
  \label{flowFrame}
\end{equation}

\end{lemma}

\noindent 
There is a rather obvious similarity between the above spacetime kinematics for the Ricci flow, the role of the lapse and shift vector field, and the use of $Dreibein$s in the formulation of the initial value problem in general relativity. However, this similarity cannot be pushed too far on the dynamical side. As a matter of fact, the \emph{natural} spacetime metric on $M^{4}_{Par}$ associated with the dynamics of the Ricci flow is not the product metric described by the diffeomorphism $F_{\eta }:I\times \Sigma\rightarrow  M^{4}_{Par}$.  Formal metrics, often strongly degenerate in the time--like direction, seem to better capture the most relevant aspects of the spacetime geometry of the Ricci flow \cite{ChowChu1,ChowChu2}, \cite{18}.

\subsection{Factorization of the linearized Ricci flow}
\label{spectra}

 \noindent As already stressed, an important role in  Ricci flow theory  is  played
by  the formal linearization of (\ref{mflow}) in the direction of a symmetric bilinear form  $h_{ab}(\beta )$,  \emph{i.e.}
\begin{equation}
\begin{tabular}{l}
$\frac{\partial  }{\partial  \beta   }{h}_{ab}(\beta  )=-\left.\frac{d}{dt}\,\left(2{\mathcal{R}_{ab}(g^{(t)})}
\right)\right|_{t=0},$ \\ 
\\ 
${h}_{ab}(\beta  =0)=h_{ab}$\, ,\;\; $0\leq \beta  <T_{0}$\;,%
\end{tabular}
  \label{linDT}
\end{equation}
where $h_{ab}(\beta )$ can be thought of as representing an infinitesimal deformation $g^{(t)}_{ab}(\beta )=
g_{ab}(\beta )+ t\;h_{ab}(\beta )$, $t\in (-\varepsilon ,\varepsilon )$, of the  flow $\beta \rightarrow g_{ab}(\beta)$ defined by (\ref{mflow}), \emph{i.e.},  $h(\beta )\,\in \mathcal{T}_{(\Sigma ,g(\beta ))}\,\mathcal{M}et(\Sigma )$. According to a lenghty but standard computation, (see \emph{e.g.}, \cite{6}, \cite{Glickenstein}, \cite{carfora:deformation1}), the linearization (\ref{linDT}) characterizes  the flow $\beta \mapsto h_{ab}(\beta )$ as a  
solution of the weakly-parabolic
initial value problem 
\begin{equation}
\begin{tabular}{l}
$\frac{\partial  }{\partial  \beta }h_{ab}=\Delta _{L}h_{ab}+2\,\left[\delta_{g} ^{*}\,\delta_{g} \,G(h)  \right]_{ab}
\;, $ \\ 
\\ 
$h_{ab}(\beta =0)=h_{ab}\;,\;\; 0\leq \beta  <T_{0}$ \;.%
\end{tabular}   \label{linflow}
\end{equation}
For notational ease, in (\ref{linflow}) we have dropped the explicit $\beta $-dependence
and we have introduced the Einstein--conjugate  $G(g,h)\doteq h-\frac{1}{2}\,\left(tr_{g}\,h \right)\,g$  of $h\in C^{\infty }(\Sigma ,\otimes ^{2}T^{*}\,\Sigma)$,
($G(h)$ for short, if it is clear, from the context, with respect to which metric $g$ we are conjugating).
\noindent The operator $\Delta _{L}\,:$ $C^{\infty }(\Sigma ,\otimes ^{2}T^{*}\,\Sigma)\rightarrow C^{\infty }(\Sigma ,\otimes ^{2}T^{*}\,\Sigma)$ is the Lichnerowicz-DeRham Laplacian on symmetric bilinear
forms defined by 
\begin{equation}
\Delta _{L}{h}_{ab}\doteq \triangle {h}%
_{ab}-R_{as}{h}_{b}^{s}-R_{bs}{h}_{a}^{s}+2R_{asbt}{h%
}^{st},
\label{LDR}
\end{equation}
where $\triangle \doteq g^{ab}(\beta )\,\nabla _{a}\,\nabla _{b}$ is the rough (or Bochner) Laplacian, and where  for $n=3$ we can set
\begin{equation}
R_{asbt}=R_{ab}g_{st}+R_{st}g_{ab}-R_{sb}g_{at}-R_{at}g_{sb}+\frac{1}{2}R\left(g_{at}g_{sb}-g_{ab}g_{st} \right)\;.
\label{treRiem}
\end{equation}
For each given $\beta\in [0,T_{0})$,  $\Delta _{L}$ is an operator of Laplace type \cite{gilkey}, \emph{i.e.}, $\Delta _{L}=\Delta+\mathcal{E}$, for $\mathcal{E}$ the local section of $End\,\left(\otimes ^{2}T^{*}\,\Sigma\right)$, \,
$h_{ik}\mapsto \mathcal{E}^{ik}_{ab}\,h_{ik}$,  provided by
\begin{equation}
\mathcal{E}^{ik}_{ab}\doteq -3\mathcal{R}_{a}^{i}\delta _{b}^{k}-3\mathcal{R}_{b}^{k}\delta _{a}^{i}+2\mathcal{R}^{ik}g_{ab}+
2\left(\mathcal{R}_{ab}-\frac{1}{2}\mathcal{R}g_{ab} \right)g^{ik}+\mathcal{R}\delta_{a}^{i}\delta_{b}^{k}\;.
\label{lichendo}
\end{equation}
$\Delta _{L}$ is $L^{2}$ self--adjoint, $\left(\Delta _{L}\,h,\,k\right)_{L^{2}(\Sigma )}$,$=\left(h,\,\Delta _{L}\,k\right)_{L^{2}(\Sigma )}$, but it is not negative semi--definite since  
\begin{eqnarray}
&&\;\;\;\;\;\int_{\Sigma }h^{ab}\,\Delta _{L}\,h_{ab}\,d\mu _{g}=\\
&&-\int_{\Sigma }\left[\nabla_{i}\,h_{ab}\,\nabla^{i}\,h^{ab}+
6\,h^{ab}\mathcal{R}_{as}h^{s}_{b}-4h\,h^{ab}\mathcal{R}_{ab}+
\mathcal{R}\left(h^{2}-h_{ab}h^{ab} \right) \right]\,d\mu _{g}\nonumber\;,
\end{eqnarray}
where $h\doteq g^{ab}\,h_{ab}$. Along the Ricci flow the curvature can grow unboundedly large, thus, in order to have some control on the spectral properties of  $\Delta _{L}$, we need to restrict attention to a particular subclass of Ricci flow metrics. In particular,
we shall say that a Ricci flow $\beta \mapsto g_{ab}(\beta )$ on $\Sigma \times [0,T_{0})$ is of bounded geometry on the subinterval $[0,\beta ^{*}]\subset [0,T_{0})$ if, in such an interval, the associated $\beta $--dependent curvature and its covariant derivatives of each order have uniform bounds, \emph{i.e.}, if there exists constants $C_{k}>0$ such that $\left|\nabla ^{k}\,Rm(\beta )\right|\leq C_{k}$, \, $k=0,1,\ldots$, \, for $0\leq \beta \leq \beta ^{*}$. The hypothesis of bounded geometry considerably simplifies the characterization of the conjugate linearized Ricci flow (in particular the analysis of the associated heat kernel and of its asymptotics), without sacrificing generality. By exploiting the technique of parabolic rescalings, one can extend the analysis to  Ricci flow singularities, at least in the case when one has a noncollapsed limit, (\emph{e.g.}, for finite time singularities on closed manifolds).

\noindent If we assume that $\beta \mapsto g_{ab}(\beta )$ on $\Sigma \times [0,T_{0})$ is of bounded geometry on the subinterval $[0,\beta ^{*}]\subset [0,T_{0})$,  then  from the spectral theory of Laplace type operators on closed Riemannian manifolds (see  \cite{gilkey}, and \cite{gilkey2} (Th. 2.3.1)), it follows that, on $(\Sigma ,g_{ab}(\beta ))$, for each given $\beta \in [0,\beta ^{*}]\subset [0,T_{0})$, the operator $P_{L}\doteq -\Delta _{L}=-\left(\Delta +\mathcal{E}\right)$  has a discrete spectral resolution $\left\{h_{ik}^{(n)}(\beta ),\,\lambda _{(n)}(\beta )\right\}$, where the  ordered eigenvalues
$\lambda _{(1)}(\beta )\leq\lambda _{(2)}(\beta )\leq \ldots\infty $ have finite multiplicities, and are contained in $[-C(\beta ),\,\infty )$ for some constant $C(\beta )$ depending from the (bounded) geometry of $(\Sigma ,g(\beta ))$. Moreover, for any $\varepsilon >0$, there exists an integer $n(\varepsilon;\beta  )$ so that $n^{\frac{2}{3}-\varepsilon }\leq \lambda _{(n)}\leq n^{\frac{2}{3}+\varepsilon }$, for $n\geq n(\varepsilon;\beta  )$.   The corresponding set of eigentensor $\left\{h_{ik}^{(n)}(\beta )\right\}$, \ $h_{ik}^{(n)}(\beta )\in C^{\infty }(\Sigma ,\otimes ^{2}T^{*}\,\Sigma)$, with $P_{L}\,h_{ik}^{(n)}(\beta )=\lambda _{(n)}(\beta )\,h_{ik}^{(n)}(\beta )$,  provide a complete orthonormal basis for $L^{2 }(\Sigma ,\otimes ^{2}T^{*}\,\Sigma)$. If for a tensor field $\phi\in L^{2 }(\Sigma ,\otimes ^{2}T^{*}\,\Sigma)$ we denote by $c_{n}(\beta )\doteq \left(\phi,h_{(n)}(\beta ) \right)_{L^{2}(\Sigma)}$ the corresponding Fourier coefficients, then  
$\phi_{ik}\,\in C^{\infty }(\Sigma ,\otimes ^{2}T^{*}\,\Sigma)$ iff $\lim_{n\rightarrow \infty }n^{k}\,c_{n}(\beta )=0$, $\forall k\in \mathbb{N}$, (\emph{i.e}, the $\{c_{n}(\beta )\}$ are rapidly decreasing). Also, if $|\phi|_{k}$ denotes the $\sup$--norm of $k^{th}$ covariant derivative of $\phi$, then there exists $j(k;\beta )$ so that $|\phi|_{k}\leq \,n^{j(k;\beta )}$ if $n$ is large enough. This result implies in particular that the series $\phi_{ab}=\sum_{n}c_{n}(\beta )\,h^{(n)}_{ab}(\beta )$ converges absolutely to $\phi_{ab}$ in the $C^{\infty }$ topology. 

\noindent In order to exploit the properties of $P _{L}$ for defining the conjugate linearized Ricci flow we need to factorize (\ref{linflow}) into a strictly parabolic flow and a $Diff(\Sigma )$ generating term. There are various distinct ways of implementing such a decomposition, all eventually related to the DeTurck trick \cite{DeTurck}. For the convenience of the reader, here  we describe a well--known factorization \cite{01} in a form particularly suited to our purposes, (to the best of my knowledge, such a factorization appeared first explicitly in \cite{lottCMC}), and which holds for any $n$--dimensional manifold. Further details can be found in (Chap.2 of) \cite{chowluni}.

\noindent  Let us consider a given symmetric bilinear form
$\widetilde{h}_{ab}^{(0)}\in T_{g}\mathcal{M}et(\Sigma )$. Along the Ricci flow of metrics $%
\beta \longmapsto g_{ab}(\beta )$, $g_{ab}(\beta=0 )=g_{ab}$, $0\leq \beta <T_{0}$, look for solutions $\beta \mapsto h_{ab}(\beta )$ of the associated linearized flow in the form 
\begin{equation}
h_{ab}(\beta )=\widetilde{h}_{ab}(\beta )+\nabla _{a}w_{b}(\beta )+\nabla _{b}w_{a}(\beta ),
\label{elledue}
\end{equation}
with $\widetilde{h}_{ab}(\beta=0 )=\widetilde{h}_{ab}^{(0)}$, and
where the $\beta $-dependent\ vector field \ $w^{a}(\beta )$ is
associated with  $\beta $-dependent infinitesimal $Diff(\Sigma )$
reparametrizations of \ the Riemannian structure associated with $%
g_{ab}(\beta )$. 

\noindent Since $\widetilde{h}_{ab}(\beta )+\mathcal{L}_{\vec{w}}g_{ab}$
must satisfy the linearized Ricci flow, we get\ 

\begin{equation}
\begin{tabular}{l}
$\frac{\partial  }{\partial  \beta }\widetilde{h}_{ab}+\frac{\partial  }{%
\partial  \beta }\mathcal{L}_{\vec{w}}g_{ab}=\Delta _{L}\widetilde{h}_{ab}-\Delta
_{L}\mathcal{L}_{\vec{w}}g_{ab} $ \\ 
\\ 
$+2\left[\delta^{\ast }_{g}\delta_{g}\,G(\widetilde{h})\right]_{ab}
+2\left[\delta^{\ast }_{g}\delta_{g}\,G(\mathcal{L}_{\vec{w}}{g})\right]_{ab}
$%
\end{tabular}
 ,  \label{explic}
\end{equation}

\noindent where $\widetilde{h}(\beta )\doteq g^{ab}(\beta )\,\widetilde{h}_{ab}(\beta )$.\ \ From the relations 
\begin{equation}
\frac{\partial  }{\partial  \beta }\mathcal{L}_{\vec{w}(\beta )}g_{ab}(\beta )=%
\mathcal{L}_{\vec{w}(\beta )}\frac{\partial  }{\partial  \beta }g_{ab}(\beta )+%
\mathcal{L}_{\frac{\partial }{\partial  \beta }\vec{w}(\beta )}g_{ab}(\beta ),
\end{equation}
\begin{equation}
\mathcal{L}_{\vec{w}(\beta )}\frac{\partial  }{\partial  \beta }g_{ab}(\beta )=-2%
\mathcal{L}_{\vec{w}(\beta )}R_{ab}
\end{equation}
(in the latter we have exploited the fact that $\beta \longmapsto
g_{ab}(\beta )$ evolves along the  Ricci flow), and 
\begin{equation}
\mathcal{L}_{\vec{w}}R_{ab}=-\frac{1}{2}\Delta _{L}\mathcal{L}_{\vec{w}}g_{ab}-\left[\delta_{g} ^{\ast }\delta_{g}\,G(\mathcal{L}_{\vec{w}}{g})\right]_{ab}\;,
\end{equation}
(consequence of the the $Diff(\Sigma )$-equivariance of the Ricci tensor) we
obtain 
\begin{equation}
\frac{\partial  }{\partial  \beta }\mathcal{L}_{\vec{w}}g_{ab}=\mathcal{L}_{\frac{%
\partial  }{\partial  \beta }\vec{w}}g_{ab}+\Delta _{L}\mathcal{L}_{\vec{w}}g_{ab}
+2\,\left[\delta_{g} ^{\ast }\delta_{g}\,G(\mathcal{L}_{\vec{w}}{g})\right]_{ab}\;. \label{liebob}\\
\end{equation}
Inserting this latter in (\ref{explic}) we have
\begin{equation} 
\frac{\partial  }{\partial  \beta }\,\widetilde{h}_{ab}+\mathcal{L}_{(\frac{%
\partial }{\partial  \beta }{w}_{k}+\nabla ^{i}(\widetilde{h}_{ik}-\frac{1}{2}%
\widetilde{h}\,g_{ik}))}\,\,g_{ab}=\Delta _{L}\widetilde{h}_{ab}\;.
\label{mainfactor}
\end{equation}
\noindent As an immediate consequence of the structure of this relation it follows that, under the stated hypotheses, we can naturally factorize the linearized Ricci flow according to the (see \emph{e.g.}, \cite{chowluni})
\begin{lemma} (The Reduced Linearized Ricci Flow).
\ Let 
$\beta \longmapsto \widetilde{h}_{ab}(\beta )$, $\beta\in [0,T _{0})$,
denote the flow solution of the parabolic initial value problem 
\begin{equation} 
\begin{tabular}{l}
$\frac{\partial  }{\partial  \beta }\,\widetilde{h}_{ab}=\Delta _{L}\widetilde{h%
}_{ab}$ \\ 
\\  
$\widetilde{h}_{ab}(\beta =0)={h}_{ab},$%
\end{tabular}
   \label{divfree}
\end{equation}
and let $\beta \longmapsto w_{a}(\beta )$, $\beta\in [0,T _{0})$, be the $%
\beta $-dependent (co)vector field solution of the  initial value problem
\begin{equation} 
\begin{tabular}{l}
$\frac{\partial  }{\partial  \beta }\,w_{a}(\beta )=-\nabla ^{b}\,\left(\widetilde{h}_{ab}-\frac{1}{2}\,\widetilde{h}\,g_{ab}\right) ,$ \\ 
\\ 
$\;\;w_{a}(\beta =0)=0,$%
\end{tabular}
   \label{tracecons}
\end{equation}
then the flow $\beta \longmapsto h_{ab}(\beta )$, $\beta\in(0,\beta _{0})$,
defined by
\begin{equation}
h_{ab}(\beta )\doteq \widetilde{h}_{ab}(\beta )+\mathcal{L}_{\vec{w}(\beta
)}g_{ab}(\beta ),
\end{equation}
solves the linearized Ricci flow (\ref{linflow}) with initial datum $%
h_{ab}(\beta =0)={h}_{ab}$. 
\label{factlemma}
\end{lemma}

\begin{proof}
The proof of the lemma amounts to backtracking the steps leading to the identity (\ref{mainfactor}). Explicitly,
from (\ref{divfree}) and (\ref{liebob}), we get

\begin{equation} 
\begin{tabular}{l}
$\frac{\partial  }{\partial  \beta }\,\widetilde{h}_{ab}+\frac{\partial }{\partial \beta }\,\mathcal{L}_{\vec{w}}\,g_{ab}=\Delta _{L}\widetilde{h}_{ab}$ \\ 
\\
$+\mathcal{L}_{\frac{%
\partial }{\partial  \beta }{w}}\,\,g_{ab}+\Delta _{L}\,\mathcal{L}_{\vec{w}}\,g_{ab}
+2\,\left[\delta_{g} ^{\ast }\delta_{g}\,G(\mathcal{L}_{\vec{w}}{g})\right]_{ab}\;. $\\
\end{tabular}
\label{mainfactor2}
\end{equation}
Moreover, from (\ref{tracecons}),  we  have 
\begin{equation}
\mathcal{L}_{\frac{%
\partial }{\partial  \beta }{w}}\,\,g_{ab}+  2\,\left[\delta_{g} ^{\ast }\delta_{g}\,G(\mathcal{L}_{\vec{w}}{g})\right]_{ab}
= 2\,\left[\delta_{g} ^{\ast }\delta_{g}\,G(\widetilde{h}+\mathcal{L}_{\vec{w}}{g})\right]_{ab}\;.
\label{mainfactor4}
\end{equation}
By inserting (\ref{mainfactor4}) in (\ref{mainfactor2}), and gathering terms, we get that 
$h_{ab}(\beta )\doteq \widetilde{h}_{ab}(\beta )+\mathcal{L}_{\vec{w}(\beta
)}g_{ab}(\beta )$, solves the linearized Ricci flow (\ref{linflow}) with initial datum $%
h_{ab}(\beta =0)={h}_{ab}$. 
\end{proof}
\noindent The net effect of curvature on the factorization of the linearized Ricci flow is most easily seen in an orthonormal frame. Since every 3--manifold is parallelizable,  we can choose orthonormal sections $\{e_{(\mu)}\}_{\mu=1,2,3}$ for $(T\Sigma, g(\beta =0)$, (locally $e_{(\mu)}|_{U}=\iota^{k}_{\mu}\partial_{i}$), such that the induced basis in $\Lambda ^{2}(T_{p}\Sigma )$, $\{e_{2}\wedge e_{3},\, e_{3}\wedge e_{1},\, e_{1}\wedge e_{2} \}$   diagonalizes the curvature tensor $Rm(g)$, \emph{i.e.},
$Rm_{2323}:=r_{1}$, $Rm_{3131}:=r_{2}$, and $Rm_{1212}:=r_{3}$. Let us evolve the sections $\{e_{(\mu)}\}_{\mu=1,2,3}$ along the given  Ricci flow according to the Uhlenbeck trick (\ref{flowHulen}) and correspondingly set
${\frak h}_{\mu\nu }(\beta ):={\widetilde h}_{jk}(\beta )\,\iota^{j}_{\mu}(\beta )\,\iota^{k}_{\nu}(\beta )$, where 
$\beta \mapsto {\widetilde h}_{jk}(\beta )$ is the solution of the reduced linearized Ricci flow (\ref{divfree}). Then, according to lemma \ref{uhlenlemma}, we get, (suppressing the $\beta$--dependence for notational ease), 
\begin{equation} 
\begin{tabular}{l}
$\frac{\partial  }{\partial  \beta }\,{\frak h}_{\mu\nu}=\Delta _{D}\,{\frak h}_{\mu\nu}+\iota^{a}_{\mu}\iota^{b}_{\nu}\,\mathcal{E}^{jk}_{ab}\,{\widetilde h}_{jk}$ \\ 
\\
$+\iota^{a}_{\mu}\iota^{b}_{\nu}\,\mathcal{R}^{c}_{a}\,{\widetilde h}_{cb}
+\iota^{a}_{\mu}\iota^{b}_{\nu}\,\mathcal{R}^{c}_{b}\,{\widetilde h}_{ac}\;. $\\
\end{tabular}
\label{}
\end{equation}
Since $\mathcal{E}^{jk}_{ab}=-\mathcal{R}_{a}^{j}\,\delta ^{k}_{b}-\mathcal{R}_{b}^{k}\,\delta ^{j}_{a}+2\mathcal{R}_{a\,b}^{j\,k}$, the above expression reduces to
\begin{equation}
\frac{\partial  }{\partial  \beta }\,{\frak h}_{\mu\nu}=\Delta _{D}\,{\frak h}_{\mu\nu}+2\,{\frak R}_{\mu\sigma\nu\tau}\,{\frak h}^{\sigma\tau}\;,
\end{equation}
where we have set ${\frak R}_{\mu\sigma\nu\tau}\doteq \iota^{a}_{\mu}\,\iota^{b}_{\nu}\,\iota^{s}_{\sigma}\,\iota^{t}_{\tau}\,\mathcal{R}_{asbt}$ and ${\frak h}^{\sigma\tau}\doteq \iota^{\sigma}_{s}\,\iota^{\tau}_{t}\,{\widetilde h}^{st}$, with $\iota^{\alpha }_{a}$  the components of the orthonormal (co)--basis $\{\theta ^{(\alpha)}\}$ dual to $\{e_{(\mu)}\}$. Thus, from Hamilton's maximum principle \cite{11a}, it follows that if $\beta \mapsto g_{ab}(\beta )$, $0\leq \beta < \beta ^{*}\subset [0,T_0)$ is a Ricci flow with non--negative curvature operator and with bounded geometry and if $\beta \mapsto {\widetilde h}_{ij}(\beta )$ is a solution of the reduced linearized Ricci flow (\ref{divfree}) with ${\widetilde h}_{ij}(\beta =0)>0$, then ${\widetilde h}_{ij}(\beta)>0$ for every $\beta \in [0,\beta ^{*}]$.

\noindent If, in the initial value problems (\ref{divfree}) and (\ref{tracecons}), we consider the  initial conditions $\widetilde{h}_{ab}(\beta =0)=0$, and $w_{a}(\beta =0)=\xi _{a}$, then one recovers the well--known fact that, for a $\beta $--independent vector  $\vec\xi\in C^{\infty }(\Sigma ,\,T\Sigma ) $, the tensor field $h_{ab}(\beta )=\mathcal{L}_{\xi }\,\,g_{ab}(\beta )$, is a solution of the linearized Ricci flow, and that any Killing vector  is preserved along the Ricci flow. More generally, the existence of the $\mathcal{D}iff\,(\Sigma )$--solitonic solutions of the Ricci flow, and the structure of the factorization described by lemma \ref{factlemma}, suggest that there should exist solutions of the reduced linearized Ricci flow (\ref{divfree}) of the form $\widetilde{h}_{ab}(\beta )=\mathcal{L}_{v(\beta ) }\,\,g_{ab}(\beta )$ for some judiciously chosen $\beta \mapsto v^{a}(\beta )$. This is expressed by the following

\begin{lemma} 
\ For a given Ricci flow $\beta \mapsto g_{ab}(\beta )$, $0\leq \beta < T_{0}$, let  $\beta \mapsto v_{a}(\beta )$ denote the flow solution of the parabolic initial value problem
\begin{equation}
\begin{tabular}{l}
$\frac{\partial  }{\partial  \beta }v_{a}(\beta )=\triangle v_{a}(\beta )+\mathcal{R}_{a}^{b}v_{b}(\beta ),$ \\ 
\\ 
$v_{a}(\beta =0)=v_{a},$%
\end{tabular}  
\label{vettv}
\end{equation}
where  ${v}(\beta=0)\,\in C^{\infty }(\Sigma ,T^{*}\Sigma )$ is a given covector field, and where  $\triangle v_{b}(\beta )+\mathcal{R}_{b}^{a}v_{a}(\beta )$, with $\triangle v_{b}(\beta )\doteq \nabla ^{a}\nabla _{a}v_{b}(\beta )$, is the (1--form) Hodge laplacian on $(\Sigma ,g(\beta ))$. Then the flow $\beta\mapsto \widetilde{h}_{ab}(\beta )=\mathcal{L}_{v(\beta ) }\,\,g_{ab}(\beta )$ provides a $\mathcal{D}iff(\Sigma )$--solitonic solution to the reduced linearized Ricci flow (\ref{divfree}).
\label{shiftlemma}
\end{lemma}

\noindent Again, in a form or another, this is a well--known property of the linearized Ricci flow, see \emph{e.g.} \cite{chowluni}, (note that in \cite{chowluni} the sign convention on Ricci tensor is opposite to ours). Here we are emphasizing, for later use, the $\mathcal{D}iff(\Sigma )$--solitonic nature of such solutions. 

\begin{proof}
A direct computation using the Ricci commutation relations provides
\begin{eqnarray}
&&-\,\delta_{g} \,G(\mathcal{L}_{\vec{v}}{g})
=\label{harm2} \\
&&=\nabla ^{a}\left[\nabla _{a}v_{b}(\beta )+\nabla _{b}v_{a}(\beta )-g_{ab}(\beta)\nabla ^{c}v_{c}(\beta )  \right]\,dx^{b}=\nonumber\\
&&=\left[\triangle v_{b}(\beta )+\mathcal{R}_{ab}v^{a}(\beta )\right]\,dx^{b}  \nonumber     \;,
\end{eqnarray}
 thus, according to (\ref{vettv})
\begin{equation}
\frac{\partial  }{\partial  \beta }\,v_{a}(\beta )=-\left[\delta_{g} \,G(\mathcal{L}_{\vec{v}}{g})\right]_{a} \;,
 \label{tracebons}
\end{equation}
which implies, (see (\ref{menoLie})),   
\begin{equation}
\mathcal{L}_{\frac{%
\partial  }{\partial  \beta }\vec{v}}g_{ab}=-2  \left[\delta_{g}^{*}\delta_{g} \,G(\mathcal{L}_{\vec{v}}{g})\right]_{ab}\;. 
\end{equation}
By introducing this latter relation in
\begin{equation}
\frac{\partial  }{\partial  \beta }\mathcal{L}_{\vec{v}}g_{ab}=\mathcal{L}_{\frac{%
\partial  }{\partial  \beta }\vec{v}}g_{ab}+\Delta _{L}\mathcal{L}_{\vec{v}}g_{ab}+
2\,\left[\delta_{g}^{*}\delta_{g} \,G(\mathcal{L}_{\vec{v}}{g})\right]_{ab}\;,
\label{liebobv}
\end{equation}
(see (\ref{liebob})), we get
\begin{equation}
\frac{\partial  }{\partial  \beta }\,\mathcal{L}_{\vec{v}}\,g_{ab}=\Delta _{L}\mathcal{L}_{\vec{v}}\,g_{ab}\;, 
\label{hlie} 
\end{equation}
\noindent which implies that  $\widetilde{h}_{ab}(\beta )=\mathcal{L}_{v(\beta )}\,g_{ab}(\beta )$ solves (\ref{divfree}) with the initial datum $\left.\mathcal{L}_{\vec{v}(\beta )}\,g_{ab}(\beta)\right|_{\beta =0}=\mathcal{L}_{\vec{v}}\,g_{ab}$.
\end{proof}

\noindent Lemma \ref{shiftlemma} and of eqn. (\ref{hlie}),\, may suggest that, along the Ricci flow, we can decompose the given solution $\beta \mapsto \widetilde{h}_{ab}(\beta )$ of (\ref{divfree}) according to 
\begin{equation}
\widetilde{h}(\beta ) = \widetilde{h}^{(T)}(\beta )+2\delta_{g} ^{*}{v(\beta )},\;\;\;\delta_{g} \,\widetilde{h}^{(T)}(\beta )=0\;.
\label{LsplitDec}
\end{equation}
This would also imply that the divergence--free part $\widetilde{h}^{(T)}(\beta )$ evolves according to $\frac{\partial }{\partial \beta }\widetilde{h}^{(T)}(\beta )=\Delta _{L}\widetilde{h}^{(T)}(\beta )$. However, from $\delta_{g} \widetilde{h}^{(T)}(\beta )=0$ it follows that the (co)vector field defined by Lemma \ref{shiftlemma} must also comply with the constraint $2\delta_{g} \delta_{g} ^{*}v=\delta_{g} \widetilde{h}(\beta )$, for all $0\leq \beta <T_{0}$, (in components this reduces to the elliptic PDE \, $\triangle v_{a}+\mathcal{R}_{ab}v^{b}+\nabla _{a}\nabla ^{b}v_{b}=\nabla ^{b}\widetilde{h}_{ab}$, \, where $\widetilde{h}_{ab}(\beta )$ is the given source). Such a requirement clearly overdetermines $\beta \mapsto v_{a}(\beta )$, and we cannot assume that (\ref{LsplitDec}) holds in the general case. This also follows explicitly from the following  
\begin{lemma} (A commutation formula) For any symmetric bilinear form $S_{kl}$ on any $n$-dimensional Riemannian manifold, we have 
\begin{eqnarray}
&&\nabla ^{k}\,\triangle _{L}\,S_{kl}=\triangle \,\nabla ^{k}\,S_{kl}+S^{ab}\,\nabla ^{k}\,R_{kalb}-R_{la}\,\nabla ^{k}\,S_{k}^{a}-S_{k}^{a}\,\nabla ^{k}\,R_{la}\label{lichcomm}\\
\notag\\
&=&\triangle \,\nabla ^{k}\,S_{kl}+S^{ka}\,\nabla _{l}\,R_{ka}-R_{la}\,\nabla ^{k}\,S_{k}^{a}-2S_{k}^{a}\,\nabla ^{k}\,R_{la}
\notag\;.
\end{eqnarray}
\end{lemma}
\begin{proof}
The proof is a somewhat lengthy but otherwise standard computation exploiting Ricci commutation formulas and the second Bianchi identity. In detail
\begin{equation*}
\begin{array}{lll}
-\nabla^k\triangle_LS_{ik} &=&\nabla^k(-\nabla^j\nabla_jS_{ik}+R^l_iS_{lk}+R^l_kS_{il}-2R^{lj}_{ik}S_{lj})\\
&=&-(\nabla^j\nabla^k(\nabla_jS_{ik})-R^k_j\nabla^jS_{ik}-R^{kl}_{ij}\nabla^jS_{kl}+R_j^k\nabla^jS_{ik})\\
&&+(\nabla^kR^l_i)S_{lk}+R^l_i\nabla^kS_{lk}+(\nabla^kR^l_k)S_{il}+R^l_k\nabla^kS_{il}-\\
&&2(\nabla^kR^{lj}_{ik})S_{lj}-2R^{lj}_{ik}\nabla^kS_{lj}\\
&=&-\nabla^j(\nabla_j\nabla^kS_{ik}-R^{kl}_{ji}S_{kl}+R_j^kS_{ik})+R^{kl}_{ij}\nabla^jS_{kl}\\
&&+(\nabla^kR^l_i)S_{lk}+R^l_i\nabla^kS_{lk}+(\nabla^kR^l_k)S_{il}+R^l_k\nabla^kS_{il}-\\
&&2(\nabla^kR^{lj}_{ik})S_{lj}-2R^{lj}_{ik}\nabla^kS_{lj}\\
&=&-\triangle\nabla^kS_{ik}-(\nabla^jR^{kl}_{ij})S_{kl}+(\nabla^kR^j_i)S_{jk}+R^j_i\nabla^kS_{jk}\;.
\end{array}
\end{equation*}
From the second Bianchi identity we get $\nabla^jR^{kl}_{ij} = -\nabla^kR^l_i+\nabla_iR^{kl}$, which inserted into the above expression eventually provides (\ref{lichcomm}).
\end{proof}
From (\ref{lichcomm}) and the Ricci flow rule
\begin{equation}
\frac{\partial }{\partial \beta }\,\nabla ^{k}\,S_{kl}= g^{ik}\,\nabla _{i}\left(\frac{\partial }{\partial \beta }\,S_{kl} \right)+    2R^{ik}\nabla _{i}\,S_{kl}+\,S^{mi}\,\nabla_{l}R_{mi}\;,
\end{equation}
(which follows directly from the evolution of the Christoffel symbols under the Ricci flow), we immediately compute that if $\beta \mapsto S_{kl}(\beta )$ evolves, along the Ricci flow, according to $\frac{\partial }{\partial \beta }\,S_{kl}(\beta )=\Delta _{L}\,S_{kl}(\beta )$, then 
\begin{eqnarray}
\label{gradS}
\frac{\partial }{\partial \beta }\,\nabla ^{k}\,{S}_{kl}&=&\triangle \,\nabla ^{k}\,{S}_{kl}
-R_{l}^{a}\,\nabla ^{k}\,{S}_{ka}+ {S}^{ab}\,\nabla _{l}\,R_{ab}+\\
\notag\\
&&+2R^{ik}\,\nabla _{i}\,{S}_{kl} - 2{S}^{ik}\,\nabla _{i}\,{R}_{kl}\;.\notag
\end{eqnarray}
The presence, in the above expression, of the terms $ {S}^{ab}\,\nabla _{l}\,R_{ab}
+2R^{ik}\,\nabla _{i}\,{S}_{kl} - 2{S}^{ik}\,\nabla _{i}\,{R}_{kl}$, 
implies that, unless we are on a ($3$--dimensional) Einstein manifold, $R_{ab}=\frac{1}{3}R\, g_{ab}$, the parabolic initial value problem (\ref{gradS}) with initial the condition $\left. \nabla ^{a}\,{S}_{ab}(\beta)\right|_{\beta =0}=0$, 
does not admit, in general, the solution $\nabla ^{a}\,{S}_{ab}(\beta)=0$,\, $0\leq \beta <T_{0}$.
If we apply this latter result to $S(\beta )=\left(2\delta_{g} ^{*}v(\beta )-\widetilde{h}(\beta )\right)$ it follows that the  $L^{2}(\Sigma _{\beta },\,g(\beta ))$--orthogonal decomposition 
\begin{equation}
C^{\infty }(\Sigma_{\beta },\,\otimes ^{2}T^{*}\Sigma_{\beta }  )\cong  Ker\;\delta_{g(\beta)} \,\oplus Im\;\delta_{g(\beta)} ^{*} \label{Lsplit}
\end{equation}
cannot be naturally imposed to the coupled evolution $\beta\mapsto (g_{ab}(\beta ), \widetilde{h}_{ab}(\beta ))$ along a generic Ricci flow metric on $\Sigma \times [0,T_{0})$.

\noindent  The difficulties one experiences in controlling the $L^{2}$--decomposition of the solutions of (\ref{divfree}) are related to the \emph{dynamical} $\mathcal{D}iff(\Sigma )$--equivariance of (\ref{linflow}) and are a counterpart of the existence of the solitonic solutions of the Ricci flow. It is then natural to bypass such difficulties  by adopting a strategy akin to the one used by G. Perelman in handling Ricci flow $\mathcal{D}iff(\Sigma )$--solitons.
In particular, in order to have an a priori control on the $L^{2}(\Sigma _{\beta },\,g(\beta ))$ decomposition (\ref{Lsplit}), we shall characterize the (backward) flow which is conjugated to the $\mathcal{D}iff(\Sigma )$--soliton solutions of (\ref{divfree}), described by lemma \ref{shiftlemma}.

\section{The conjugate linearized Ricci flow}
\label{LRF}
Let $\beta \mapsto (\Sigma ,\,g_{ab}(\beta ))$, $0\leq \beta \leq \beta ^{*}$, $\beta ^{*}\in [0,T_{0})$ be a given Ricci flow metric of bounded geometry, and let $(M^{4}_{Par}\simeq \Sigma \times [0,\beta ^{*}],\,{}^{(4)}g_{par})$ denote the corresponding parabolic spacetime. Through the diffeomorphism $F^{-1}_{\beta }:M^{4}_{par} \rightarrow I\times \Sigma $, (see (\ref{Fdiff})), any $(\beta,x) \mapsto B_{ab}(\beta,x )$, with 
$B_{ab}(\beta,x )\in C^{\infty }(\Sigma,\,\otimes ^{2}T^{*}\Sigma)$, can be seen as an element of the space of symmetric bilinear forms on $M^{4}_{Par}$, $C^{\infty }(M^{4}_{Par},\,\otimes ^{2}T^{*}M^{4}_{Par})$.
Since the volume form on $M^{4}_{Par}$ is given by the product measure $d\mu _{g(\beta )}\,d\beta $, we can consider, on 
$C^{\infty }(M^{4}_{Par},\,\otimes ^{2}T^{*}M^{4}_{Par})$,\,
the $L^{2}(M^{4}_{Par},\,{}^{(4)}g_{par})$ inner product
\begin{equation}
\int_{0}^{\beta ^{*}}\,\int_{\Sigma }\,g^{ia}(\beta) g^{kb}(\beta )\,H_{ik}(\beta )\,B_{ab}(\beta )\,d\mu _{g(\beta )}\,d\beta \;,
\end{equation}
between $H_{ik}(\beta )$ and  $B_{ab}(\beta )$ $\in C^{\infty }(\Sigma,\,\otimes ^{2}T^{*}\Sigma)$. Similarly, we can also define the natural pairing 
\begin{equation}
\int_{0}^{\beta ^{*}}\,\int_{\Sigma }\,H^{ab}(\beta )\,B_{ab}(\beta )\,d\mu _{g(\beta )}\,d\beta \;,
\end{equation}
between $H^{ab}(\beta )\in  C^{\infty }(\Sigma ,\,\otimes ^{2}T\Sigma )$ and $B_{ab}(\beta )\in C^{\infty }(\Sigma ,\,\otimes ^{2}T^{*}\Sigma )$.\, Let us consider the operator
\begin{equation} 
\bigcirc _{L}  \doteq \frac{\partial }{\partial \beta }-\triangle _{L}\;,
\end{equation}
acting on the space of $\beta $--dependent symmetric bilinear forms $C^{\infty }(\Sigma ,\,\otimes ^{2}T^{*}\Sigma )$\, $\subset $\, $C^{\infty }(M^{4}_{Par},\,\otimes ^{2}T^{*}M^{4}_{Par})$. According to lemma \ref{shiftlemma},  $Ker\,\,\bigcirc _{L}\cap \,Im\,\delta_{g} ^{*}$ characterizes the solitonic solutions of the (reduced) linearized Ricci flow (\ref{divfree}).
Let us compute its 
$L^{2}(M^{4}_{Par},\,{}^{(4)}g_{par})$ conjugate $\bigcirc^{*} _{L}$, thought of as acting on the space  of symmetric two--tensors with compact support.  From the relation
\begin{eqnarray}
&&\int_{0}^{\beta ^{*}}\,\int_{\Sigma }\,H^{ab}(\beta )\,\frac{\partial }{\partial \beta }B_{ab}(\beta )\,d\mu _{g(\beta )}\,d\beta =\int_{0}^{\beta ^{*}}\,\frac{d}{d\beta }\int_{\Sigma }\,H^{ab}\,B_{ab}\,d\mu _{g}\,d\beta\\
\notag\\
&&+\int_{0}^{\beta ^{*}}\,\int_{\Sigma }\,B_{ab} \left(-\frac{\partial }{\partial \beta }H^{ab}+\mathcal{R}\,H^{ab}\right) \,d\mu _{g}\,d\beta \notag\\
\notag\\
&&=\int_{0}^{\beta ^{*}}\,\int_{\Sigma }\,B_{ab} \left(-\frac{\partial }{\partial \beta }H^{ab}+\mathcal{R}\,H^{ab}\right)\,d\mu _{g}\,d\beta \notag\;,
\end{eqnarray}
(where we have exploited the Ricci flow evolution for $d\mu _{g}$ and the time--boundary condition $H^{ab}\in C_{0}^{\infty }(M^{4}_{Par},\,\otimes ^{2}TM^{4}_{Par})$), and 
\begin{eqnarray}
&&\int_{0}^{\beta ^{*}}\,\int_{\Sigma }\,H^{ab}(\beta )\,\left(-\triangle _{L} \right)\,B_{ab}(\beta )\,d\mu _{g(\beta )}\,d\beta \\
\notag\\
&&=\int_{0}^{\beta ^{*}}\,\int_{\Sigma }\,B_{ab}(\beta )\,\left(-\triangle _{L} \right)\,H^{ab}(\beta )\,d\mu _{g(\beta )}\,d\beta \notag\;,
\end{eqnarray}
(where we have exploited the fact that $\triangle _{L}$ is formally self--adjoint on each $(\Sigma ,g(\beta ))$), we  compute
\begin{eqnarray}
&&\int_{0}^{\beta ^{*}}\,\int_{\Sigma }\,H^{ab}(\beta )\,\bigcirc _{L}\,B_{ab}(\beta )\,d\mu _{g(\beta )}\,d\beta \\
\notag\\
&&=\int_{0}^{\beta ^{*}}\,\int_{\Sigma }\,H^{ab}\,\left(\frac{\partial }{\partial \beta }-\triangle _{L}\ \right)\,B_{ab}\,d\mu _{g}\,d\beta \notag\\
\notag\\
&&=\int_{0}^{\beta ^{*}}\,\int_{\Sigma }\,B_{ab}\,\left(-\frac{\partial }{\partial \beta }-\triangle _{L} +\mathcal{R} \right)\,H^{ab}\,d\mu _{g}\,d\beta\notag\\
\notag\\
&&=\int_{0}^{\beta ^{*}}\,\int_{\Sigma }\,B_{ab}(\beta )\,\bigcirc^{*} _{L}\,H^{ab}(\beta )\,d\mu _{g(\beta )}\,d\beta\notag\;.
\end{eqnarray}
Thus,
\begin{equation}
\bigcirc^{*} _{L}\doteq -\frac{\partial }{\partial \beta }-\triangle _{L} +\mathcal{R}\;.
\label{conjflow}
\end{equation}

\noindent The following results provide the geometrical meaning of $\bigcirc^{*} _{L}$.
\begin{lemma}
Let $\beta \mapsto (\Sigma ,g(\beta ))$, $0\leq \beta \leq \beta ^{*}$, be a  Ricci flow of bounded geometry, and 
let $Ker\;\delta_{g} $ denote the corresponding $\beta $--dependent subspace of divergence--free $2$--tensor fields $H^{ab}(\beta )\in C^{\infty }(\Sigma ,\,\otimes ^{2}T\Sigma )$,  then $Ker\,\delta_{g} $ is an invariant subspace for $\bigcirc^{*} _{L}$, \emph{i.e.}
\begin{equation}
\bigcirc^{*} _{L}\left(Ker\,\delta_{g}  \right)\subset Ker\,\delta_{g} \;,
\end{equation}  
for all $\beta \in [0,\beta ^{*}]$.
\label{includiv}
\end{lemma}

\begin{proof}
The commutation formula (\ref{lichcomm}) provides
\begin{equation}
\nabla _{a}\triangle _{L}H^{ab}=\triangle \nabla _{a}H^{ab}+H^{aj}\nabla ^{b}\mathcal{R}_{aj}-
\mathcal{R}^{b}_{a}\nabla _{j}H^{aj}-2H^{aj}\nabla _{j}\mathcal{R}^{b}_{a}\;,
\end{equation}
whereas along the Ricci flow we have 
\begin{eqnarray}
&&-\nabla _{a}\frac{\partial }{\partial \beta }H^{ab}=-\frac{\partial }{\partial \beta }\left(\nabla _{a}H^{ab}\right)\\
\notag\\
&&-H^{rb}\nabla _{r}\mathcal{R}-H^{ar}\nabla _{a}\mathcal{R}^{b}_{r}-H^{ar}\nabla _{r}\mathcal{R}^{b}_{a}
+H^{ra}\nabla ^{b}\mathcal{R}_{ar}\;\notag.
\end{eqnarray}
Inserting these relations into the expression for $\nabla _{a}\left(\bigcirc^{*} _{L}\,H^{ab}\right)$, and cancelling terms, we easily get
\begin{eqnarray}
\label{divcirc}
&&\nabla _{a}\left(\bigcirc^{*} _{L}\,H^{ab}\right)=-\nabla _{a}\left[\left(\frac{\partial }{\partial \beta }+\triangle _{L} -\mathcal{R}\right)H^{ab}\right]\\
\notag\\
&&=-\left(\frac{\partial }{\partial \beta }+\triangle -\mathcal{R}\right)\nabla _{a}H^{ab}+\mathcal{R}^{b}_{a}\nabla _{j}H^{aj}\;,\notag
\end{eqnarray}
(note that the Laplacian in the last line is the rough Laplacian). Thus, if $\nabla _{a}H^{ab}(\beta )=0$ then 
$\nabla_{a}\left(\bigcirc^{*} _{L}\,H^{ab}(\beta )\right)=0$.
\end{proof}

\noindent As expected under $L^{2}$--duality,  the action of $\bigcirc ^{*}_{L}$ on $Im\, \delta_{g} ^{*}$ parallels the rather complicate action of $\bigcirc _{L}$ on $Ker\, \delta_{g}$. In particular,  for the Lie derivative  along a gradient vector field $X^{a}(\eta )\doteq g^{ak}\,\nabla _{k}\,f$, with $f\in C^{\infty }(\Sigma ,\,\mathbb{R}) $, we have 

\begin{lemma}
Let $\left(Hess\,f(\beta )\right)^{ab}\doteq g^{ia}(\beta )g^{kb}(\beta )\,\nabla _{i}\nabla _{k}\,f(x,\beta )$, be the (contravariant) Hessian of a $\beta $--dependent function $f\in C^{\infty }(\Sigma ,\,\mathbb{R}) $, then along the Ricci flow $\beta \mapsto g_{ab}(\beta )$ on $\Sigma \times [0,\beta ^{*}]$ we have

\begin{eqnarray}
&&\bigcirc _{L}^{*}\, \left(Hess\,f(\beta )\right)^{ab}=\,\nabla ^{a}\,\nabla ^{b}\left(\frac{\partial }{\partial \beta }+\Delta -\mathcal{R} \right)\,f+\label{hess}\\
\nonumber\\
&&+2\,\nabla _{i}\left(\mathcal{R}^{ia}\,\nabla ^{b}f \right)+2\,\nabla _{k}\left(\mathcal{R}^{kb}\,\nabla ^{a}f \right)+\nonumber\\
\nonumber\\
&&+2\,\left(\nabla ^{a}\mathcal{R}^{b}_{l}+\nabla ^{b}\mathcal{R}^{a}_{l}-\nabla _{l}\mathcal{R}^{ab} \right)\,\nabla ^{l}f\nonumber\;.
\end{eqnarray}
\end{lemma}
\noindent Note in particular that the component of $\bigcirc _{L}^{*}\, \left(Hess\,f(\beta )\right)$ on $Im\, \delta_{g} ^{*}$ is generated by the $L^{2}(M^{4}_{Par})$--adjoint
\begin{equation}
\square^{*}\doteq -\left(\frac{\partial}{\partial\beta} +\Delta -\mathcal{R}\right)\;, 
\label{adheatop}
\end{equation}
of the scalar heat operator 
\begin{equation}
\square\doteq \left(\frac{\partial}{\partial\beta} -\Delta \right)\;.
\label{heatop} 
\end{equation}

\begin{proof} The proof of (\ref{hess}) is a long but routine computation exploiting the Ricci flow identity
\begin{eqnarray}
&&\nabla _{i}\,\nabla _{k}\left(\frac{\partial }{\partial \beta }+\Delta \right)\,f=
\left(\frac{\partial }{\partial \beta }+\Delta_{L} \right)\,\nabla _{i}\,\nabla _{k}f \\
\nonumber\\
&&-2\,\left(\nabla _{i}\mathcal{R}_{kl}+\nabla _{k}\mathcal{R}_{il}-\nabla _{l}\mathcal{R}_{ik} \right)\,\nabla ^{l}f\nonumber\;,
\end{eqnarray}
(for this latter see \cite{chowluni}, Chap. 2, \, \S 5).
\end{proof}
\noindent Consider the set of covector fields $\vec{v}(\beta )\in C^{\infty }(M^{4}_{Par},\,T^{*}M^{4}_{Par})$ obtained as solutions of 
\begin{equation} 
\begin{tabular}{l}
$\frac{\partial  }{\partial  \beta }v_{a}(\beta )=\triangle v_{a}(\beta )+\mathcal{R}_{a}^{b}v_{b}(\beta ),$ \\ 
\\ 
$\;\;v_{a}(\beta =0)=v^{(0)}_{a},$%
\end{tabular}  
\label{vettv2}
\end{equation}
where the initial  ${\vec{v}}_{(0)}$ varies in $C^{\infty }(\Sigma ,T^{*}\Sigma )$. According to lemma \ref{shiftlemma}, these flows describe all possible solitonic solutions   $\mathcal{L}_{\vec{v}(\beta )}g_{ab}(\beta )$ of the linearized Ricci flow (\ref{divfree}). Let ${H}^{ab}(\beta )$, $\beta\in [0,\beta ^{*}]$, be a $\beta $--dependent  2--tensor field, and  let us consider the pairing  
\begin{equation}
\int_{\Sigma }{H}^{ab}(\beta )\,\mathcal{L}_{\vec{v}(\beta )}\,g_{ab}(\beta )\,d\mu _{g(\beta )}\;,
\label{Tortho}
\end{equation}
for every $0\leq \beta \leq \beta ^{*}$. By differentiating (\ref{Tortho}), and exploiting (\ref{hlie} ),  we get
\begin{eqnarray}
&&\frac{d}{d\beta }\int_{\Sigma }{H}^{ab}(\beta )\,\mathcal{L}_{\vec{v}(\beta )}\,g_{ab}(\beta )\,d\mu _{g(\beta )}\\
\notag\\
&&=\int_{\Sigma }\left[\mathcal{L}_{\vec{v}}\,g_{ab}\,\frac{\partial }{\partial \beta }{H}^{ab}\,+ {H}^{ab}\left(\triangle _{L} -\mathcal{R} \right)\,\mathcal{L}_{\vec{v}}\,g_{ab}\right]d\mu _{g}\notag\\
\notag\\
&&=-\int_{\Sigma }\mathcal{L}_{\vec{v}}\,g_{ab}\,\bigcirc ^{*}_{L}\,{H}^{ab}\,d\mu _{g}\notag\;,
\end{eqnarray}
which implies
\begin{eqnarray}
&&\left.\int_{\Sigma }{H}^{ab}(\beta )\,\mathcal{L}_{\vec{v}(\beta )}\,g_{ab}(\beta )\,d\mu _{g(\beta )}\right|_{\beta ^{*}}-
\left.\int_{\Sigma }{H}^{ab}(\beta )\,\mathcal{L}_{\vec{v}(\beta )}\,g_{ab}(\beta )\,d\mu _{g(\beta )}\right|_{\beta=0}\\
\notag\\
&&=-\int_{0}^{\beta ^{*}}\int_{\Sigma }\mathcal{L}_{\vec{v}}\,g_{ab}\,\bigcirc ^{*}_{L}\,{H}^{ab}\,d\mu _{g}\,d\beta \notag\;.
\end{eqnarray}

\noindent Thus, if we evolve ${H}^{ab}(\beta )$   according to the flow 
\begin{equation}
\bigcirc ^{*}_{L}\,{H}^{ab}(\beta )=0\;,
\label{circevol}
\end{equation}
the inner product (\ref{Tortho}) will be preserved along the evolution, \emph{i.e.}
\begin{equation}
\left.\int_{\Sigma }{H}^{ab}\,\mathcal{L}_{\vec{v}}\,g_{ab}\,d\mu _{g}\right|_{\beta ^{*}}=
\left.\int_{\Sigma }{H}^{ab}\,\mathcal{L}_{\vec{v}}\,g_{ab}\,d\mu _{g}\right|_{\beta=0}\;.
\end{equation}
Since any solution $\beta \mapsto g_{ab}(\beta )$ of the Ricci flow on $\Sigma_{\beta } \times [0,\beta ^{*}]$  can be converted into a solution $\eta \mapsto g_{ab}(\eta )$ of the backward Ricci flow on $\Sigma_{\eta } \times [0,\beta ^{*}]$  by the time reversal $\beta \mapsto \eta \doteq \beta ^{*}-\beta $, the above remarks motivate the following
\begin{definition}
Let $\eta \mapsto g_{ab}(\eta )$, \ $\eta \doteq \beta ^{*}-\beta $, a backward Ricci flow on $\Sigma_{\eta } \times [0,\beta ^{*}]$, then the  conjugated evolution $\bigcirc ^{*}_{L}\,{H}^{ab}=0$, of a symmetric bilinear form $H^{ab}(\eta =0)$, along $\eta \mapsto g_{ab}(\eta )$ takes the form of the parabolic initial value problem  
\begin{equation} 
\begin{tabular}{l}
$\frac{\partial }{\partial \eta }{H}^{ab}=\Delta_{L}{H}^{ab}-\,\mathcal{R}{H}^{ab}\;,$\\
\\
${H}^{ab}(\eta=0)={H}^{ab}_{*}\;.$%
\end{tabular}
\;   \label{transVol}
\end{equation}  
\end{definition}
\noindent Note that, according to the backward $\beta $--parabolic character of the operator $\bigcirc ^{*}_{L}(\beta )$, the initial data ${H}^{ab}(\eta=0)={H}^{ab}_{*}$  in (\ref{transVol}) correspond to $\beta =\beta ^{*}$.
Lemma \ref{includiv} trivially extends to the evolution (\ref{transVol}) and we have the
\begin{corollary}
\label{Hdivf}
If $\eta \longmapsto H^{ab}(\eta )$, $0\leq \eta  \leq \beta ^{*}$, is the solution of
the parabolic initial value problem (\ref{transVol}) with $\nabla _{a}\,{H}^{ab}_{*}=0$, then $\nabla _{a}\,{H}^{ab}(\eta)=0$, $\forall \eta \in [0,\beta ^{*}]$.
\end{corollary}

\noindent Moreover, we have the following result that explicitly shows that (\ref{transVol}) is conjugated to (\ref{divfree}).

\begin{proposition}
Let $\eta \mapsto {{H}}^{ab}(\eta )$, $\eta\in [0,\beta ^{*}]$,\;${H}^{ab}(\eta =0)$ $={H}^{ab}_{*}$, be a solution of the parabolic initial value problem (\ref{transVol}). Also, let 
$\beta \mapsto {\widetilde{h}}_{ab}(\beta )$, $\beta\in [0,\beta ^{*}]$,\,${\widetilde{h}}_{ab}(\beta =0)=h_{ab}(\beta=0)$ be a solution of reduced linearized Ricci flow (\ref{divfree}). Then, along the backward Ricci flow 
$\eta \mapsto g_{ab}(\eta )$, \ $\eta \doteq \beta ^{*}-\beta $, on $\Sigma_{\eta } \times [0,\beta ^{*}]$,
the  flows $\eta \mapsto {{H}}^{ab}(\eta)$ and $\eta \mapsto {\widetilde{h}}_{ab}(\eta ):={\widetilde{h}}_{ab}(\beta^{*}-\beta )$ are $L^{2}(M^{4}_{Par})$ conjugated, in the sense that
\begin{equation}
\frac{d}{d\eta }\int_{\Sigma }{H}^{ab}(\eta )\,\widetilde{h}_{ab}(\eta )\,d\mu _{g(\eta )}=0\;.
\end{equation}
In particular,  
\begin{equation}
\left.\int_{\Sigma }{H}^{ab}(\eta )\,\widetilde{h}_{ab}(\eta )\,d\mu _{g(\eta )}\right|_{\eta =0}=
\left.\int_{\Sigma }{H}^{ab}(\eta )\,\widetilde{h}_{ab}(\eta )\,d\mu _{g(\eta )}\right|_{\eta =\beta ^{*}}\;.
\end{equation}
\label{usprop}
\end{proposition}

\begin{proof}
A direct computation provides
\begin{eqnarray}
&&\;\;\;\;\frac{d}{d\eta }\int_{\Sigma }{H}^{ab}(\eta )\,\widetilde{h}_{ab}(\eta )\,d\mu _{g(\eta )}
\label{longrel}\\
\notag\\
&=&\int_{\Sigma }\left(\Delta_{L}{H}^{ab}-\,\mathcal{R}{H}^{ab} \right)\widetilde{h}_{ab}\,d\mu _{g}
+\int_{\Sigma }H^{ab}\left(-\Delta _{L}\widetilde{h%
}_{ab}+\widetilde{h}_{ab}\,\mathcal{R} \right)\,\,d\mu _{g}=0
\;.\notag
\end{eqnarray}
\end{proof}

\noindent This result and corollary \ref{Hdivf} directly imply the

\begin{proposition} 
let $\eta \mapsto (\Sigma ,g(\eta ))$,\, $0\leq \eta \leq \beta ^{*}$, be a backward Ricci flow of bounded geometry. Assume that $C^{\infty }(\Sigma_{\eta },\,\otimes ^{2}T^{*}\Sigma_{\eta }  )\supset  Ker\,\delta_{g} \not=\emptyset$, \, $0\leq \eta \leq \beta ^{*}$. Let $\eta \mapsto {{H}}^{ab}_{(T)}(\eta )$, \, $\nabla _{a}{{H}}^{ab}_{(T)}(\eta )=0$,\,   $\eta\in [0,\beta ^{*}]$,\;${H}^{ab}_{(T)}(\eta =0)$ $={H}^{ab}_{*}$, with $\delta_{g}\,{H}^{ab}_{*}=0$, be a divergence--free solution of the parabolic initial value problem (\ref{transVol}). If 
$\beta \mapsto {\widetilde{h}}_{ab}(\beta )$, $\beta\in [0,\beta ^{*}]$,\,${\widetilde{h}}_{ab}(\beta =0)=h_{ab}(\beta=0)$ denotes a solution of reduced linearized Ricci flow (\ref{divfree}), and $\eta \mapsto \widetilde{h}_{ab}^{(T)}(\eta ):=\widetilde{h}_{ab}^{(T)}(\beta^{*}-\beta )$ is its divergence--free part along $\eta \mapsto (\Sigma ,g(\eta ))$,  then
\begin{equation}
\int_{\Sigma }{H}^{ab}_{(T)}(\eta)\,\widetilde{h}_{ab}^{(T)}(\eta )\,d\mu _{g(\eta )}\;,
\end{equation}
is constant along the coupled backward evolution $\eta \mapsto \left(g_{ab}(\eta ),\widetilde{h}_{ab}(\eta )\right)$.
\end{proposition}

\begin{proof}
By writing $\widetilde{h}_{ab}(\eta )=\widetilde{h}_{ab}^{(T)}(\eta )+\mathcal{L}_{X(\eta )}\,g_{ab}(\eta )$, for some $\eta $--dependent vector field $X(\eta )$, and exploiting the $L^{2}$--orthogonality between ${H}^{ab}_{(T)}(\eta)$ and  $\mathcal{L}_{X(\eta )}\,g_{ab}(\eta )$, we get
$\int_{\Sigma }{H}^{ab}_{(T)}(\eta)\,\widetilde{h}_{ab}(\eta )\,d\mu _{g(\eta )}$\,$=\int_{\Sigma }{H}^{ab}_{(T)}(\eta)\,\widetilde{h}_{ab}^{(T)}(\eta )\,d\mu _{g(\eta )}$.
\end{proof}

\noindent Thus, the conjugated flow (\ref{transVol}) provides the \emph{directions} in  $C^{\infty }(\Sigma_{\beta },\,\otimes ^{2}T^{*}\Sigma_{\beta } )$\; along which the non-trivial solutions $\beta \mapsto \widetilde{h}_{ab}^{(T)}(\beta )$ of the linearized Ricci flow (\ref{divfree}) propagate without dissipation in the $L^{2}$ sense. In this connection notice also that along the the conjugated flow (\ref{transVol}) we have the following monotonicity result

\begin{proposition}
 Let $\eta \mapsto g_{ab}(\eta )$, \ $\eta \doteq \beta ^{*}-\beta $, a backward Ricci flow of bounded geometry on $\Sigma_{\eta } \times [0,\beta ^{*}]$ with $\mathcal{R}(\eta )\geq 0$, $\eta \doteq \beta ^{*}-\beta $, where $\mathcal{R}(\eta)$ denotes the scalar curvature of $(\Sigma, g(\eta) )$.   If $\eta \mapsto {{H}}^{ab}(\eta )$, $\eta\in [0,\beta ^{*}]$,\;${H}^{ab}(\eta =0)$ $={H}^{ab}_{*}$, denotes a solution of the parabolic initial value problem (\ref{transVol}) with $\delta _{g(\eta )}H(\eta)\not=0$, then
 \begin{equation}
 \frac{d}{d\eta }\,\int_{\Sigma }\left| \delta _{g(\eta )}H(\eta)\right|^{2}\,d\mu _{g(\eta )}\leq 0\;,
 \end{equation}
 where $\left| \delta _{g(\eta )}H(\eta)\right|^{2}\doteq \nabla _{a}{H}^{ab}(\eta)\,\nabla _{c}{H}^{cd}(\eta )\,g_{bd}(\eta)$.
 \end{proposition}
\begin{proof}
From (\ref{divcirc}) we get
\begin{equation}
\label{vecconj}
\frac{\partial }{\partial \eta }\nabla _{a}H^{ab}=\Delta \nabla _{a}H^{ab} -\mathcal{R}\,\nabla _{a}H^{ab}-\mathcal{R}^{b}_{a}\nabla _{j}H^{aj}\;,
\end{equation}
from which we compute
\begin{eqnarray}
&&\frac{\partial }{\partial \eta }\,\left| \delta _{g(\eta )}H(\eta)\right|^{2}=2g_{bd}\,\nabla _{c}H^{cd}\,\Delta\left(\nabla _{a}H^{ab} \right) -2 \mathcal{R}\,\left| \delta _{g(\eta )}H(\eta)\right|^{2}\\
&&=\Delta\,\left| \delta _{g(\eta )}H(\eta)\right|^{2}-2\, \left|\nabla\,\delta _{g(\eta )}H(\eta) \right|^{2} 
 -2 \mathcal{R}\,\left| \delta _{g(\eta )}H(\eta)\right|^{2}\nonumber    \;,
\end{eqnarray}
where $\left|\nabla\,\delta _{g(\eta )}H(\eta) \right|^{2}\doteq \nabla _{i} \left(\nabla _{a}{H}^{ab}\right)\nabla ^{i}\,\left(\nabla _{c}{H}^{cd}\right)\,g_{bd}$. By integrating, and taking into account that along the backward Ricci flow  $\frac{\partial }{\partial \eta }\,d\mu _{g(\eta )}=\mathcal{R}\,d\mu _{g(\eta )}$, we get

\begin{eqnarray}
\frac{d}{d\eta }\,\int_{\Sigma }\left| \delta _{g(\eta )}H(\eta)\right|^{2}\,d\mu _{g(\eta )}=
-\int_{\Sigma }\mathcal{R}\,\left| \delta _{g(\eta )}H(\eta)\right|^{2}\,d\mu _{g(\eta )}\\
\nonumber\\
-2\, \int_{\Sigma }\left|\nabla\,\delta _{g(\eta )}H(\eta) \right|^{2}\,d\mu _{g}
\leq 0\;.\nonumber
\end{eqnarray}
\end{proof}
\noindent Since non--negative scalar curvature is preserved along the Ricci flow, the requirement $\mathcal{R}(\eta )\geq 0$, $\eta \doteq \beta ^{*}-\beta $, in the above result is not particularly restrictive. In particular, it can be easily removed by weighting the riemannian measure $d\mu _{g}$ with a positive solution of the forward conjugate scalar heat equation, (I wish to thank Lei Ni for this latter remark). According to proposition \ref{usprop}, it also follows that  (\ref{vecconj}) is the backward  flow $L^{2}(M^{4}_{Par})$--conjugated to the forward evolution for covector fields defined by lemma \ref{shiftlemma}. 

\noindent These elementary aspects of the $L^{2}(M^{4}_{Par})$ conjugacy relation have an important and rather unexpected consequence, which implies that the conjugate linearized Ricci flow averages out the full Ricci flow:  
\begin{proposition} 
Let $\beta \mapsto (\Sigma ,{g}(\beta ))$, $\beta\in [0,\beta ^{*}]$ be a  Ricci flow of bounded geometry, and 
let $\beta \mapsto {\mathcal{R}}_{ab}(\beta )$ be the corresponding $\beta $--evolution of the Ricci tensor. Denote by
$\eta \mapsto {{H}}^{ab}(\eta )$, $\eta\in [0,\beta ^{*}]$, \;${H}^{ab}(\eta =0)$ $={H}^{ab}_{*}$ the solution of the parabolic initial value problem (\ref{transVol}) associated with the given $\beta \mapsto (\Sigma ,{g}(\beta ))$. Then, 
\begin{equation}
\label{bello1}
\frac{d}{d\eta }\int_{\Sigma }{H}^{ab}(\eta)\,\mathcal{R}_{ab}(\eta )\,d\mu _{g(\eta )}=0\;,
\end{equation}
and 
\begin{equation}
\label{bello2}
\frac{d}{d\eta }\int_{\Sigma }\left({g}_{ab}(\eta )-2\eta\,\mathcal{R}_{ab}(\eta ) \right){H}^{ab}(\eta)\,d\mu _{g(\eta )}=0\;,
\end{equation}
along the backward Ricci flow.
In particular, this implies 
\begin{equation}
\left.\int_{\Sigma }{H}^{ab}\,\mathcal{R}_{ab}\,d\mu _{g}\right|_{\eta =0}
=\int_{\Sigma }{H}^{ab}(\eta)\,\mathcal{R}_{ab}(\eta )\,d\mu _{g(\eta )}\;,
\end{equation}
and 
\begin{equation}
\left.\int_{\Sigma }{H}^{ab}\,{g}_{ab}\,d\mu _{g}\right|_{\eta =0}
=\int_{\Sigma }\left({g}_{ab}(\eta )-2\eta\,\mathcal{R}_{ab}(\eta ) \right){H}^{ab}(\eta)\,d\mu _{g(\eta )}\;,
\label{usprop2a}
\end{equation}
for every $0\leq \eta \leq \beta ^{*}$.
\label{usprop2}
\end{proposition}

\begin{proof}
It is easily checked that in any dimension $n$  the forward evolution for the Ricci curvature
\begin{equation}
\frac{\partial }{\partial \beta }{\mathcal{R}}_{ij}=\Delta
{\mathcal{R}}_{ij}-6g^{kl}{%
\mathcal{R}}_{il}{\mathcal{R}}_{kj}+3\mathcal{R}{\mathcal{R}}%
_{ij}
+2g_{ij}{\mathcal{R}}^{kl}{%
\mathcal{R}}_{kl}-g_{ij}\,\mathcal{R}^{2}\;, \label{scalev3}
\end{equation}
can be expressed directly in terms of the Lichnerowicz--DeRham Laplacian as 
\begin{equation}
\frac{\partial }{\partial \beta }{\mathcal{R}}_{ij}=\Delta_{L}
{\mathcal{R}}_{ij}\;. \label{scalev4}
\end{equation}
Thus, from (\ref{transVol}) we get
\begin{eqnarray}
&&\frac{d}{d\eta }\int_{\Sigma }{H}^{ab}(\eta )\,\mathcal{R}_{ab}(\eta )\,d\mu _{g(\eta )}\\
&&=\int_{\Sigma }\left\{\mathcal{R}_{ab}\frac{\partial }{\partial \eta }H^{ab}+H^{ab}\frac{\partial }{\partial \eta }\mathcal{R}_{ab}+ {H}^{ab}\,\mathcal{R}_{ab}\,\mathcal{R}    \right\}  \notag\\
&&=\int_{\Sigma }\left\{\mathcal{R}_{ab}\left[\triangle _{L}H^{ab}
-\mathcal{R}\,{H}^{ab} \right]   \right. \notag\\
&&-\left.\,H^{ab}\,\triangle _{L}\mathcal{R}_{ab}+{H}^{ab}\,\mathcal{R}_{ab}\,\mathcal{R} \right\}\,d\mu _{g}\notag\\
&&=\int_{\Sigma }\left\{-H^{ab}\,\triangle _{L}\,\mathcal{R}_{ab}+\mathcal{R}_{ab}\,\triangle _{L}\,H^{ab}\right\}\,d\mu _{g}=0\;,\notag
\end{eqnarray}
from which (\ref{bello1}) follows. Relation (\ref{bello2}) follows similarly by observing that, since $\eta\mapsto g_{ab}(\eta)$ is covariantly constant, we can write
\begin{equation}
\frac{\partial }{\partial\eta }\,\left({g}_{ab}(\eta )-2\eta\,\mathcal{R}_{ab}(\eta ) \right)=-\,\Delta _{L}\,\left({g}_{ab}(\eta )-2\eta\,\mathcal{R}_{ab}(\eta ) \right)\;.
\end{equation} 
\end{proof}
\noindent This result has an interesting converse
\begin{remark}
Let $\eta \mapsto g_{ab}(\eta )$\;$\in \mathcal{M}et(\Sigma )$ be a one--parameter family of evolving metrics on $\Sigma \times [0,\beta ^{*}]$, not identified a priori with a backward Ricci flow. Let $\eta \mapsto \left(g_{ab}(\eta )\,,H^{ab}(\eta)\right)$ be the corresponding solution of the heat equation 
\begin{equation}
\begin{tabular}{l}
$\eta \mapsto g_{ab}(\eta )\;,$\;\;\;$0\leq \eta \leq \beta ^{*}$\;,\\
\\
$\frac{\partial }{\partial \eta }{H}^{ab}=\Delta_{L}{H}^{ab}-\,\mathcal{R}{H}^{ab}\;,$\\
\\
${H}^{ab}(\eta=0)={H}^{ab}_{*}\;,$\;\;\; ${H}^{ab}_{*}\in C^{\infty }(\Sigma ,\otimes ^{2}\,T\Sigma )$\;.
\end{tabular}
\;   \label{L2formul}
\end{equation} 
Then among all possible such flows $\eta \mapsto \left(g_{ab}(\eta )\,,H^{ab}(\eta)\right)$, the backward Ricci flow
$\eta \mapsto g_{ab}(\eta )$,\,$\frac{\partial }{\partial \eta }\,g_{ab}(\eta )=2\,\mathcal{R}_{ab}$ is characterized by
the  condition
\begin{equation}
\frac{d}{d\eta }\,\int_{\Sigma }\,g_{ab}H^{ab}\,d\mu _{g(\eta )}= \,2\,\int_{\Sigma }\,\mathcal{R}_{ab}H^{ab}\,d\mu _{g(\eta )}\;.
\end{equation}

\end{remark}

\begin{proof}
A direct computation provides
\begin{eqnarray}
&&\frac{d}{d\eta }\,\int_{\Sigma }\,g_{ab}H^{ab}\,d\mu _{g(\eta )}\\
\nonumber\\
&&=\int_{\Sigma }\,H^{ab}\left(\frac{\partial }{\partial \eta }
g_{ab}-\mathcal{R}g_{ab}+\frac{1}{2}g_{ab}g^{ik}\frac{\partial }{\partial \eta }
g_{ik} \right)\,d\mu _{g(\eta )}\nonumber\;,
\end{eqnarray}
which yields $2\,\int_{\Sigma }\,\mathcal{R}_{ab}H^{ab}\,d\mu _{g(\eta )}$, for every possible solution $\eta \mapsto H^{ab}(\eta )$ of (\ref{L2formul}), iff
\begin{equation}
\frac{\partial }{\partial \eta }g_{ab}(\eta )=2\,\mathcal{R}_{ab}\;.
\end{equation}

\end{proof}

It is important to stress that actually the above results (as most results in this paper) hold in any dimension $n\geq3$, this true in particular for Propositions 3.5, 3.6, 3.7, and 3.8.\\

\noindent Now we turn to the analysis of the  conjugate flow $\eta \mapsto (g_{ab}(\eta),H^{ab}(\eta ))$ in its role as the Ricci flow integral kernel. 

\section{The conjugate backward heat kernel}

The averaging properties of the conjugate linearized Ricci flow become manifest when we identify the flow $\eta \mapsto H^{ab}(\eta )$ with the heat kernel of $\bigcirc ^{*}_{L}$ along the backward Ricci flow $\eta \mapsto g_{ab}(\eta )$. To fix notation, let $(\otimes ^{2}T\Sigma)\,\boxtimes\,( \otimes ^{2}T^{*}\Sigma)$ denote the bundle over $\Sigma \times \Sigma $ whose fiber over $(y,x)$ $\in \Sigma \times \Sigma $ is given by $\left(\otimes ^{2}T\Sigma\,\boxtimes\, \otimes ^{2}T^{*}\Sigma\right)_{(y,x)}$\,$=$\,
$(\otimes ^{2}T\Sigma)_{y}\,\otimes \,( \otimes ^{2}T^{*}\Sigma)_{x}$. Wheras for notational simplicity we keep on assuming $n=3$,
it is perhaps appropriate to stress here once more that the results which follow actually hold in any dimension $n\geq3$, with the obvious changes in the range of tensorial indices involved.
Let $U_{\beta }\subset (\Sigma_{\beta},g(\beta ))$ be a geodesically convex neighborhood 
containing the generic point $x\in\Sigma_{\beta} $. For a chosen base point $y\in U_{\beta }$, let $l_{\beta }(y,x)$ denote the unique 
$g(\beta )$--geodesic segment $x=\exp_{y}\,u$,\;with $u\in T_{y}\Sigma $, connecting $y$ to $x$. Parallel transport along $l_{\beta }(y,x)$ 
allows to define a canonical isomorphism between the tangent space 
$T_{y}\Sigma_{\beta } $ and $T_{x}\Sigma_{\beta } $ which maps any given vector 
$\vec{v}(y)\in T_{y}\Sigma_{\beta }$ into a corresponding vector 
$\vec{v}_{P_{l_{\beta }(y,x)}}\in T_{x}\Sigma_{\beta }$. If $\{e_{(h)}(x) \}_{h=1,2,3}$ 
and $\{e_{(k')}(y) \}_{k'=1,2,3}$ respectively denote basis vectors in 
$T_{x}\Sigma_{\beta } $ and $T_{y}\Sigma_{\beta } $, (henceforth, primed indexes will always refer to components of elements of the tensorial algebra over $T_{y}\Sigma_{\beta } $), then the components of 
$\vec{v}_{P_{l_{\beta }(y,x)}}$ can be expressed as 
\begin{equation}
\left(v_{P_{l_{\beta }(y,x)}}\right)^{k}(x)= \tau ^{k}_{h'}(y,x;\beta )\,v^{h'}(y)\;,
\end{equation}
where $\tau ^{k}_{h'}$ $\in T\Sigma_{\beta }\boxtimes T^{*}\Sigma_{\beta } $  denotes the bitensor 
 associated with the parallel transport 
along $l_{\beta }(y,x)$. The Dirac $p$--tensorial measure  in 
$U_{\beta }\subset (\Sigma_{\beta},g(\beta ))$ is defined according to
\begin{equation}
\delta ^{k_{1}\ldots k_{p}}_{h_{1}'\ldots h_{p}'}(y,x;\beta ):= \otimes _{(\alpha =1)}^{p}\,\tau ^{k_{\alpha }}_{h'_{\alpha }}(y,x;\beta )\,\,\delta_{\beta } (y,x)\; ,
\end{equation}
where $\delta_{\beta } (y,x)$ is the standard Dirac measure over the Riemannian manifold 
$(\Sigma_{\beta } ,g({\beta }))$ (see \cite{lichnerowicz}). If $(\Sigma, g_{ab}(\eta ))$ is a smooth solution to the backward Ricci flow on $\Sigma_{\eta } \times [0,\beta ^{*}]$ with bounded curvature, then we can consider the $g(\eta )$--dependent  fundamental solution ${K}^{ab}_{i'k'}(y,x;\eta )$ to the conjugate heat equation (\ref{transVol}), \emph{i.e.},
\begin{equation} 
\begin{tabular}{l}
$\left(\,\frac{\partial }{\partial \eta }-\,\Delta_{L}^{(x)}+\,\mathcal{R}\,\right)\;{K}^{ab}_{i'k'}(y,x;\eta )=0\;,$\\
\\
$\lim_{\;\eta \searrow 0^{+}}\;{K}^{ab}_{i'k'}(y,x;\eta )={\delta }^{ab}_{i'k'}(y,x;)\;,$%
\end{tabular}
\;   \label{fund}
\end{equation} 
where $(y,x;\eta )\in (\Sigma \times \Sigma \backslash Diag(\Sigma \times \Sigma ))\times [0,\beta ^{*}]$, $\eta \doteq \beta ^{*}-\beta $,  $\Delta_{L}^{(x)}$ denotes the  Lichnerowicz--DeRham laplacian with respect to the variable $x$, and ${K}^{ab}_{i'k'}(y,x;\eta )$\,  is a smooth section of $(\otimes ^{2}T\Sigma)\boxtimes (\otimes ^{2}T^{*}\Sigma)$. The Dirac initial condition is understood in the distributional sense, \emph{i.e.}, for any smooth symmetric bilinear form with compact support $w^{i'k'}\in C^{\infty }_{0}(\Sigma ,\otimes ^{2}T\Sigma )$, 
\begin{equation}
\int_{\Sigma _{\eta }}{K}^{ab}_{i'k'}(y,x;\eta )\;w^{i'k'}(y)\;d\mu^{(y)} _{g(\eta )}\rightarrow w^{ab}(x)\;\;\;as\;\;\eta \searrow 0^{+}\;, 
\end{equation}
where the limit is meant in the uniform norm on  $C^{\infty }_{0}(\Sigma ,\otimes ^{2}T\Sigma )$. 
\noindent Note that the elliptic generator, associated with $\frac{\partial }{\partial \eta }-\,\Delta_{L}^{(x)}+\,\mathcal{R}$, is the operator of Laplace type on $(\Sigma ,g(\eta ))$ defined by $\Delta_{\eta } +\mathcal{F}(\eta )$\,$=\Delta _{L}-\mathcal{R}$ where $\Delta_{\eta }$ is the rough Laplacian on $(\Sigma ,g(\eta ))$. The $\eta $--dependent  endomorphism $\mathcal{F}(\eta ):C^{\infty }(\Sigma_{\eta } ,\otimes ^{2}T^{*}\,\Sigma_{\eta })\rightarrow C^{\infty }(\Sigma_{\eta } ,\otimes ^{2}T^{*}\,\Sigma_{\eta })$ is related to the endomorphism $\mathcal{E}$,  characterizing $\Delta _{L}$,  by $\mathcal{F}^{ik}_{ab}=\mathcal{E}^{ik}_{ab}-\mathcal{R}\,\delta ^{ik}_{ab}$, \emph{i.e.}, (see (\ref{lichendo})), 
\begin{equation}
\mathcal{F}^{ik}_{ab}(\eta )\doteq -3\mathcal{R}_{a}^{i}\delta _{b}^{k}-3\mathcal{R}_{b}^{k}\delta _{a}^{i}+2\mathcal{R}^{ik}g_{ab}+
2\left(\mathcal{R}_{ab}-\frac{1}{2}\mathcal{R}g_{ab} \right)g^{ik}\;,
\end{equation}
where all geometric quantities refer to $(\Sigma ,g(\eta ))$. In analogy with the spectral properties of the Lichnerowicz--DeRham Laplacian recalled in \S \ref{spectra},  the spectral theorem \cite{gilkey} implies that the operator
\begin{equation}
P_{\eta }\doteq -\left(\Delta_{\eta } +\mathcal{F}(\eta ) \right)=-\Delta _{L}+\mathcal{R}\;,
\end{equation}
has, for each given  $\eta \in [0,\beta ^{*}]$,\, 
 a discrete, finite multiplicity, spectral resolution $\left\{\phi ^{ik}_{(n)}(\eta ),\,\lambda _{(n)}(\eta )\right\}$, with 
$\lambda _{(1)}(\eta )\leq\lambda _{(2)}(\eta )\leq \ldots\infty $ contained in $[-C(\eta ),\,\infty )$, where the constant $C(\eta )$ depends from the geometry of $(\Sigma ,g(\eta ))$, and where
$\left\{\phi ^{ik}_{(n)}(\eta )\right\}$, $\phi ^{ik}_{(n)}(\eta )\in C^{\infty }(\Sigma_{\eta } ,\otimes ^{2}T^{*}\,\Sigma_{\eta })$, with
\begin{equation}
-\left(\Delta_{\eta } +\mathcal{F}(\eta ) \right)\,\phi ^{ik}_{(n)}(\eta )=\,\lambda _{(n)}(\eta )\,\phi ^{ik}_{(n)}(\eta )\;,
\end{equation}
denotes the corresponding set of eigentensors providing a complete orthonormal basis for $L^{2 }(\Sigma_{\eta } ,\otimes ^{2}T^{*}\,\Sigma_{\eta })$. The $\eta$--dependence of $\left\{\phi ^{ik}_{(n)}(\eta ),\,\lambda _{(n)}(\eta )\right\}$ makes the characterization of ${K}^{ab}_{i'k'}(y,x;\eta )$ via the spectral theorem (see \emph{e.g.}, \cite{gilkey} and \cite{Amman}) very delicate, and to prove the existence of ${K}^{ab}_{i'k'}(y,x;\eta )$ is preferable to exploit parametrix--deformation methods. These are readily available since,  along a backward Ricci flow on $\Sigma_{\eta }\times [0,\beta ^{*}]$ with bounded geometry, the metrics $g_{ab}(\eta )$ are uniformly bounded above and below for $0\leq \eta \leq \beta ^{*}$, and it does not really matter which metric we use in topologizing the spaces $C^{\infty }(\Sigma_{\eta } ,\otimes ^{2}T^{*}\,\Sigma_{\eta })$. In particular,
heat kernels  for generalized Laplacians, such as $\Delta_{\eta } +\mathcal{F}(\eta )$, (smoothly) depending on a one--parameter family of metrics $\varepsilon \mapsto g_{ab}(\varepsilon )$, $\varepsilon \geq 0$, are briefly dealt with in \cite{vergne}. The delicate setting where the parameter dependence is, as in our case, identified with the parabolic time driving the diffusion of the kernel, is discussed  in \cite{guenther2}, \cite{{chowluni}}, (see Appendix A, \S 7 for a characterization of the parametrix of the heat kernel in such a case), and in a remarkable paper by N. Garofalo and E. Lanconelli \cite{lanconelli}. Strictly speaking, in all these works, the analysis is confined to the scalar laplacian, possibly with a potential term, but the theory readily extends to generalized laplacians, always under the assumption that the metric $g_{ab}(\beta )$ is smooth as $\nearrow \beta ^{*}$. In particular, the case of generalized Laplacian on vector bundles with time--varying geometry has been studied in considerable detail by P. Gilkey and collaborators \cite{gilkey3}, \cite{gilkey4}. By adapting to our more general setting the methods used in \cite{guenther2} and in  \cite{chowluni}, when treating the scalar time-dependent Laplacian, we get
the following
\begin{theorem}
\label{theo1}
Along a backward Ricci flow on $\Sigma_{\eta }\times [0,\beta ^{*}]$ with bounded geometry, there exists a unique fundamental solution $\eta\longmapsto {K}^{ab}_{i'k'}(y,x;\eta )$ of the tensorial heat operator 
$\left(\,\frac{\partial }{\partial \eta }-\,\Delta_{L}^{(x)}+\,\mathcal{R}\,\right)$.   
\end{theorem}
\begin{proof}
The proof, (kindly provided by Stefano Romano), is a quite lengthy construction of the the heat kernel of a time-dependent generalized Laplacian. It is presented in the appendix.
\end{proof} 
\noindent  The kernel ${K}^{ab}_{i'k'}(y,x;\eta )$ is singular as $\eta_{0} \rightarrow 0$, and
the general strategy for discussing its $\eta \searrow 0^{+}$ asymptotics is, again, to model the corresponding parametrix around the Euclidean heat kernel $\left(4\pi \,\eta \right)^{-\frac{3}{2}}\,\exp\left(-\frac{d^{2}_{0}(y,x)}{4\eta } \right)$
defined in $T_{y}\Sigma$ by means of the exponential mapping associated with the initial manifold $(\Sigma ,g_{ab}(\eta =0))$. To this end,
denote by $d_{\eta }(y,x)$ the (locally Lipschitz)  distance function on $(\Sigma ,g_{ab}(\eta ))$ and by $inj\,(\Sigma ,g(\eta ))$ the associated injectivity radius. Adopt, with respect to the metric $g_{ab}(\eta )$, geodesic polar cordinates about $y\in \Sigma $, \emph{i.e.}, $x^{j'}=d_{\eta }(y,x)\,u^{j'}$, with $u^{j'}$ coordinates on the unit sphere $\mathbb{S}^{2}\subset T_{y}\Sigma $.
By adapting the analysis in \cite{{chowluni}}, \cite{lanconelli}, and \cite{gilkey3}, \cite{gilkey4} to (\ref{fund}) we have that, as $\eta \searrow 0^{+}$, and for all $(y,x)\in \Sigma $ such that $d_{0}(y,x)<\, inj\,(\Sigma ,g(0))$, there exists a sequence of smooth sections ${\Phi [h]\, }^{ab}_{i'k'}\,(y,x;\eta )$\, $\in C^{\infty }(\Sigma\times \Sigma ' ,\otimes ^{2}T\Sigma\boxtimes \otimes ^{2}T^{*}\Sigma)$,   with ${\Phi [0]\, }^{ab}_{i'k'}\,(y,x;\eta )={\tau }^{ab}_{i'k'}\,(y,x;\eta )$, such that
\begin{equation}
\frac{\exp\left(-\frac{d^{2}_{0}(y,x)}{4\eta } \right)}{\left(4\pi \,\eta \right)^{\frac{3}{2}}}\;\sum_{h=0}^{N}\eta ^{h}{\Phi [h]\, }^{ab}_{i'k'}\,(y,x;\eta )\;,
\label{uniasi}
\end{equation}
\noindent is uniformly asymptotic to ${K}^{ab}_{i'k'}(y,x;\eta )$, \;\, \emph{i.e.}, 
\begin{eqnarray}
 &&\left|{K}^{ab}_{i'k'}(y,x;\eta )- \frac{\exp\left(-\frac{d^{2}_{0}(y,x)}{4\eta } \right)}{\left(4\pi \,\eta \right)^{\frac{3}{2}}}\;\sum_{h=0}^{N}\eta ^{h}{\Phi [h]\, }^{ab}_{i'k'}\,(y,x;\eta )\right|_{\eta \searrow 0^{+}}\\
&&  =O\left(\eta ^{N-\frac{1}{2}} \right)\;,\nonumber
\end{eqnarray}
in the uniform norm on $C^{\infty }(\Sigma\times \Sigma ' ,\otimes ^{2}T\Sigma\times \otimes ^{2}T^{*}\Sigma)$.
A detailed presentation of the $\eta \searrow 0^{+}$ asymptotics of generalized Laplacians on vector bundles with time--varying geometries is discussed in \cite{gilkey3}, \cite{gilkey4}. It is worthwhile recalling that the asymptotics for the Laplace Beltrami operator plays a key role in discussing Li--Yau--Hamilton type inequalities for the scalar conjugate heat equation in Ricci flow theory (see e.g. \cite{ecker}, \cite{LeiNi}, \cite{18}).

\noindent The heat kernel ${K}^{ab}_{i'k'}(y,x;\eta )$ can be naturally normalized along the $\eta$--expanding soliton on $\mathbb{S}^{3}$ according to
\begin{lemma}
Let $\bar{g}_{ab}$ the round metric on the unit $3$--sphere $\mathbb{S}^{3}$, and, for $\eta\in[0,\beta^{*}]$, let $\eta \mapsto 4\,(T_{0}-\beta^{*}+\eta)\,\bar{g}_{ab}$ be the expanding Ricci soliton on $\mathbb{S}^{3}$ with initial radius $r(\eta=0)=2\,\sqrt{T_{0}-\beta^{*}}$ and final radius $r(\eta=\beta^{*})=2\,\sqrt{T_{0}}$. Then, along such a backward Ricci flow the heat kernel ${K}^{ab}_{i'k'}(y,x;\eta )$ scales according to
\begin{equation}
\frac{r(\eta)^{3}}{3}\,\int_{\Sigma}{\bar g}^{i'k'}(y)\,{K}^{ab}_{i'k'}(y,x;\eta )\,{\bar g}_{ab}(x)\,d{\bar\mu}_{g(x)}=1\;,
\end{equation}
where $d{\bar\mu}_{g(x,\eta)}$ is the volume element on $(\mathbb{S}^{3},\bar{g}_{ab})$.
\end{lemma}

\begin{proof}
From proposition \ref{usprop2} we get that
along the backward Ricci flow $\eta \mapsto 4\,(T_{0}-\beta^{*}+\eta)\,\bar{g}_{ab}$ we can write, for all $0\leq \eta \leq \beta ^{*}$,
\begin{eqnarray}
\mathcal{R}_{i'k'}(y,\eta=0 )&=&\lim_{\,\eta \nearrow 0^{+}}\int_{\Sigma }{K}^{ab}_{i'k'}(y,x;\eta )\,\mathcal{R}_{ab}(x,\eta )\,d\mu _{g(\eta )}\label{propvaf}\\
\nonumber\\
&=&\int_{\Sigma }{K}^{ab}_{i'k'}(y,x;\eta )\,\mathcal{R}_{ab}(x,\eta )\,d\mu _{g(\eta )}\;.\nonumber
\end{eqnarray}
Since the Ricci tensor is scale invariant we have $\mathcal{R}_{i'k'}(y,\eta=0 )=2\,\bar{g}_{i'k'}(y,\eta=0)$ and 
$\mathcal{R}_{ab}(x,\eta)=2\,\bar{g}_{ab}(x,\eta)$, moreover $d\mu _{g(x,\eta)}=r^{3}(\eta)\,d{\bar\mu} _{g(x)}$. By inserting these expressions in (\ref{propvaf}), and tracing both members with respect to $\bar{g}^{i'k'}(y,\eta=0)$, we get the stated result. 
\end{proof}
\noindent Under natural assumptions on the curvature of the supporting backward Ricci flow, the kernel ${K}^{ab}_{i'k'}(y,x;\eta )$ also exhibits point--wise positivity properties according to
\begin{lemma}
If $(\Sigma, g_{ab}(\eta ))$ is a smooth solution to a backward Ricci flow of bounded geometry on $\Sigma_{\eta } \times [0,\beta ^{*}]$ with non--negative curvature operator, then ${K}^{ab}_{i'k'}(y,x;\eta )$, $0\leq  \eta \leq \beta^{*}$, is a positive integral kernel, i.e., ${K}^{ab}_{i'k'}(y,x;\eta )v^{i'}(y)v^{k'}(y)$, \,$\forall v\in\,T_y\Sigma$, is a positive--definite quadratic form at $T^*_{(x,\eta)}\Sigma$, \,for any $(x,\eta)\in\Sigma\times[0,\beta^*]$.
\end{lemma}

\begin{proof}
We exploit the Uhlenbeck trick  in order to rewrite the evolution for ${K}^{ab}_{i'k'}(y,x;\eta )$ in a form making the proof of the positivity of ${K}^{ab}_{i'k'}(y,x;\eta )$ manifest under the stated assumptions. To this end,
choose orthonormal sections $\{e_{(\mu)}\}_{\mu=1,2,3}$ for $(T\Sigma, g(\eta =0)$, (locally $e_{(\mu)}|_{U}=\iota^{k}_{\mu}\partial_{i}$), and let us denote by $\iota^{\alpha }_{a}$  the components of the orthonormal (co)--basis $\{\theta ^{(\alpha)}\}$  dual to $\{e_{(\mu)}\}$. It is easily seen that the evolution along the backward time $\eta$ of $\iota^{k}_{\mu}$ and $\iota^{\alpha }_{a}$, consistent with the forward $\beta$ evolution (\ref{flowHulen}) of an orthonormal basis, is provided by  
\begin{equation}
\frac{\partial}{\partial\eta}\,\iota^{k}_{\mu}=-\mathcal{R}^{k}_{h}\,\iota^{h}_{\mu}\;,\;\;\;\;\;
\frac{\partial}{\partial\eta}\,\iota^{\alpha }_{a}=\mathcal{R}^{k}_{a}\,\iota^{\alpha }_{k}\;.
\end{equation}
With these preliminary remarks along the way, let us define
\begin{equation}
\mathbb{K}^{\alpha\beta}_{\gamma'\delta'}(y,x;\eta)\doteq \iota^{\alpha }_{a}(x,\eta)\,\iota^{\beta }_{b}(x,\eta)\,{K}^{ab}_{c'd'}(y,x;\eta )\,\iota^{c'}_{\gamma'}(y,\eta=0)\,\iota^{d'}_{\delta'}(y,\eta=0)\;,
\end{equation}
and consider the $\eta$--evolution of  $\mathbb{K}^{\alpha\beta}_{\gamma'\delta'}(y,x;\eta)$, (note that the primed indices do not carry $\eta$--dependence since the orthonormal basis vectors $\{\iota^{c'}_{\gamma'}\,\partial_{c'}\}$ refer to the fixed spacetime point $(y,\eta=0)$).
From the defining equation (\ref{fund}) and lemma \ref{uhlenlemma}, (applied to the $\eta$--evolution), we get, (suppressing the $\eta$--dependence for notational ease), 

\begin{equation} 
\begin{tabular}{l}
$\frac{\partial  }{\partial  \eta }\,\mathbb{K}^{\alpha\beta}_{\gamma'\delta'}
=\Delta _{D}\,\mathbb{K}^{\alpha\beta}_{\gamma'\delta'}-\mathcal{R}\,\mathbb{K}^{\alpha\beta}_{\gamma'\delta'}$ \\ 
\\
$+\,\mathcal{E}^{\alpha\beta}_{\gamma\delta}\,\mathbb{K}^{\gamma\delta}_{\gamma'\delta'}+\iota^{a}_{\mu}\,\mathcal{R}^{k}_{a}\,\iota^{\alpha}_{k}
\,\mathbb{K}^{\mu\beta}_{\gamma'\delta'}
+\iota^{a}_{\mu}\,\mathcal{R}^{k}_{a}\,\iota^{\beta}_{k}\,\mathbb{K}^{\alpha\mu}_{\gamma'\delta'} $\;,
\end{tabular}
\label{}
\end{equation}
with $\mathcal{E}^{\alpha\beta}_{\gamma\delta}(x,\eta)\doteq \iota^{\alpha }_{a}(x,\eta)\,\iota^{\beta }_{b}(x,\eta)\,\mathcal{E}^{ab}_{cd}(x,\eta)\,\iota^{c}_{\gamma}(x,\eta)\,\iota^{d}_{\delta}(x,\eta)$. 
Since $\mathcal{E}^{jk}_{ab}=-\mathcal{R}_{a}^{j}\,\delta ^{k}_{b}-\mathcal{R}_{b}^{k}\,\delta ^{j}_{a}+2\mathcal{R}_{a\,b}^{j\,k}$, the above expression reduces to
\begin{equation}
\frac{\partial  }{\partial  \eta }\,\mathbb{K}^{\alpha\beta}_{\gamma'\delta'}=
\Delta _{D}\,\mathbb{K}^{\alpha\beta}_{\gamma'\delta'}
+2\,\mathfrak{R}^{\alpha\beta}_{\gamma\delta}\,\mathbb{K}^{\gamma\delta}_{\gamma'\delta'}-\mathcal{R}\,\mathbb{K}^{\alpha\beta}_{\gamma'\delta'}\;,
\end{equation}
where we have set $\mathfrak{R}^{\alpha\beta}_{\gamma\delta}(x,\eta)\doteq \iota^{\alpha }_{a}(x,\eta)\,\iota^{\beta }_{b}(x,\eta)\,\mathcal{R}^{ab}_{cd}(x,\eta)\,\iota^{c}_{\gamma}(x,\eta)\,\iota^{d}_{\delta}(x,\eta)$.
For $\eta\searrow 0^{+}$, $\mathbb{K}^{\alpha\beta}_{\gamma'\delta'}$ approaches, in the distributional sense, the positive integral kernel ${\delta}^{\alpha\beta}_{\gamma'\delta'}(x,y;\eta=0)$, thus, Hamilton's maximum principle \cite{11a} implies that if $\eta \mapsto g_{ab}(\eta )$, $0\leq \eta < \beta ^{*}\subset [0,T_0)$, is a backward Ricci flow  with non--negative curvature operator, then $\mathbb{K}^{\alpha\beta}_{\gamma'\delta'}$, and consequently ${K}^{ab}_{i'k'}(y,x;\eta )$, remains a positive integral kernel for every $\eta \in (0,\beta ^{*}]$. 
\end{proof}

\section{An Integral representation of the Ricci flow}

\noindent We are now in position to apply proposition \ref{usprop2} to the heat kernel solution of (\ref{fund}). We have
\begin{proposition}
Let $\eta \mapsto g_{ab}(\eta )$ be a backward Ricci flow with bounded geometry on $\Sigma_{\eta }\times [0,\beta ^{*}]$, and let ${K}^{ab}_{i'k'}(y,x;\eta )$ be the (backward) heat kernel of the corresponding conjugate linearized Ricci operator 
$\bigcirc ^{*}_{L}\, {K}^{ab}_{i'k'}(y,x;\eta )=0$, for  $\eta \in (0,\beta ^{*}]$, \, 
with ${K}^{ab}_{i'k'}(y,x;\eta\searrow 0^{+} )={\delta }^{ab}_{i'k'}(y,x)$. Then 
\begin{equation}
\mathcal{R}_{i'k'}(y,\eta=0 )=\int_{\Sigma }{K}^{ab}_{i'k'}(y,x;\eta )\,\mathcal{R}_{ab}(x,\eta )\,d\mu _{g(x,\eta )}
\;,
\label{riccikernel}
\end{equation}
for all $0\leq \eta \leq \beta ^{*}$.
Moreover, as $\eta\searrow 0^{+} $, we have the uniform asymptotic expansion
\begin{eqnarray}
&&\;\;\;\;\;\mathcal{R}_{i'k'}(y,\eta=0 )=\label{ricciasi}\\
\nonumber\\
&&\frac{1}{\left(4\pi \,\eta \right)^{\frac{3}{2}}}\,\int_{\Sigma }\exp\left(-\frac{d^{2}_{0}(y,x)}{4\eta } \right)\,{\tau }^{ab}_{i'k'}(y,x;\eta )\,\mathcal{R}_{ab}(x,\eta )\,d\mu _{g(x,\eta )}\nonumber\\
\nonumber\\
&&+\sum_{h=1}^{N}\frac{\eta ^{h}}{\left(4\pi \,\eta \right)^{\frac{3}{2}}}\,\int_{\Sigma }\exp\left(-\frac{d^{2}_{0}(y,x)}{4\eta } \right)\,{\Phi [h] }^{ab}_{i'k'}(y,x;\eta )\,\mathcal{R}_{ab}(x,\eta )\,d\mu _{g(x,\eta )}\nonumber\\
\nonumber\\
&&+O\left(\eta ^{N-\frac{1}{2}} \right)\nonumber\;,
\end{eqnarray}

\noindent where ${\tau }^{ab}_{i'k'}(y,x;\eta )$ $\in T\Sigma_{\eta }\boxtimes T^{*}\Sigma_{\eta } $  is the parallel transport operator
associated with  $(\Sigma ,g(\eta ))$.
\label{prob1}
\end{proposition}

\begin{proof}
From proposition \ref{usprop2} we get that
along the backward Ricci flow on
$\Sigma \times [0,\beta ^{*}]$,  we can write, for all $0\leq \eta \leq \beta ^{*}$,
\begin{eqnarray}
\mathcal{R}_{i'k'}(y,\eta=0 )&=&\lim_{\,\eta \nearrow 0^{+}}\int_{\Sigma }{K}^{ab}_{i'k'}(y,x;\eta )\,\mathcal{R}_{ab}(x,\eta )\,d\mu _{g(\eta )}\label{avprop}\\
\nonumber\\
&=&\int_{\Sigma }{K}^{ab}_{i'k'}(y,x;\eta )\,\mathcal{R}_{ab}(x,\eta )\,d\mu _{g(\eta )}\;.\nonumber
\end{eqnarray}
Since the asymptotics (\ref{uniasi}) is uniform, we can integrate term by term, and by isolating the lower order term, we immediately get (\ref{ricciasi}). 
\end{proof}

\noindent This results illustrates the averaging properties of the backward conjugated heat kernel $\eta \mapsto {K}^{ab}_{i'k'}(y,x;\eta )$  for the Ricci curvature of the forward  flow $\beta \mapsto g_{ab}(\beta )$.
More explicitly, since $\mathcal{R}_{ab}(x,\eta )=\mathcal{R}_{ab}(x,\beta^{*}-\eta )$ and $d\mu _{g(\eta )}=d\mu _{g({\beta^{*}}-\,\eta)}$, we can equivalently rewrite (\ref{riccikernel}) along the forward Ricci flow  as
\begin{equation}
\mathcal{R}_{i'k'}(y,\beta ^{*})=\left.\int_{\Sigma }{K}^{ab}_{i'k'}(y,x;\,(\beta^{*}-\beta))\,\mathcal{R}_{ab}(x,\beta )\,d\mu _{g(\beta )}\right|_{\beta =0}\;,
\label{ricrep}
\end{equation}
which expresses the Ricci tensor at the point $y$ and at time $\beta =\beta ^{*}$ as a backward heat kernel average of the initial Ricci tensor. 

\noindent Note that a representation structurally similar to (\ref{ricciasi}) holds also for the solution $\widetilde{h}_{i'k'}(y, \eta =0)$ of the linearized Ricci flow (\ref{divfree}), \emph{i.e.},

\begin{eqnarray}
&&\;\;\;\;\;\widetilde{h}_{i'k'}(y,\eta=0 )=\label{accaasi}\\
\nonumber\\
&&\frac{1}{\left(4\pi \,\eta \right)^{\frac{3}{2}}}\,\int_{\Sigma }\exp\left(-\frac{d^{2}_{0}(y,x)}{4\eta } \right)\,{\tau }^{ab}_{i'k'}(y,x;\eta )\,\widetilde{h}_{ab}(x,\eta )\,d\mu _{g(x,\eta )}\nonumber\\
\nonumber\\
&&+\sum_{h=1}^{N}\frac{\eta ^{h}}{\left(4\pi \,\eta \right)^{\frac{3}{2}}}\,\int_{\Sigma }\exp\left(-\frac{d^{2}_{0}(y,x)}{4\eta } \right)\,{\Phi [h] }^{ab}_{i'k'}(y,x;\eta )\,\widetilde{h}_{ab}(x,\eta )\,d\mu _{g(x,\eta )}\nonumber\\
\nonumber\\
&&+O\left(\eta ^{N-\frac{1}{2}} \right)\nonumber\;.
\end{eqnarray}

\noindent By exploiting again proposition \ref{usprop2} it is also straightforward to provide an integral representation of the full
Ricci flow  in terms of the heat kernel ${K}^{ab}_{i'k'}(y,x;\,\eta)$.  Since 
$\lim_{\,\eta \searrow 0^{+}}\int_{\Sigma }{K}^{ab}_{i'k'}(y,x;\eta )\,{g}_{ab}(x,\eta )\,d\mu _{g(\eta )}$\,$=$\, ${g}_{i'k'}(y,\eta=0)$, the identity 
(\ref{usprop2a}), applied to ${K}^{ab}_{i'k'}(y,x;\,\eta))$, directly provides the

\begin{proposition}
Under the same hypotheses of proposition \ref{prob1} we have the following integral representation of the backward Ricci flow on $\Sigma _{\eta }\times (0,\beta ^{*}]$

\begin{equation}
 g_{i'k'}\,(y,\eta=0)=\int_{\Sigma }{K}^{ab}_{i'k'}(y,x;\eta  )\,\left[{g}_{ab}(x,\eta)-2\eta\,\,\mathcal{R}_{ab}(x,\eta)\right]\,d\mu _{g(x,\eta )}\;,
\label{grepres}
\end{equation}
for all $0\leq\eta\leq\beta^{*}$.

\noindent Moreover, as $\eta\searrow 0^{+}$, we have the  asymptotics
\begin{eqnarray}
&&\;\;\;\;\;\;\;\;{g}_{i'k'}(y,\eta=0)=\label{giasi}\\
\nonumber\\
&&=\frac{1}{\left(4\pi \,\eta \right)^{\frac{3}{2}}}\,\int_{\Sigma }e^{-\frac{d^{2}_{0}(y,x)}{4\eta }}\,{\tau }^{ab}_{i'k'}(y,x;\eta )\,\left[{g}_{ab}(x,\eta)-2\eta\,\mathcal{R}_{ab}(x,\eta)\right]\,d\mu _{g(x,\eta)}     \nonumber\\
\nonumber\\
&&+\sum_{h=1}^{N}\frac{(\eta) ^{h}}{\left(4\pi \,\eta \right)^{\frac{3}{2}}}\,\int_{\Sigma }e^{-\frac{d^{2}_{0}(y,x)}{4\eta }}\,{\Phi [h] }^{ab}_{i'k'}(y,x;\eta)\left[{g}_{ab}(x,\eta)-2\eta\,\mathcal{R}_{ab}(x,\eta)\right]\,d\mu _{g(x,\eta)} \nonumber\\
\nonumber\\
&&+O\left((\eta) ^{N-\frac{1}{2}} \right)\nonumber\;,
\end{eqnarray}
where $d^{2}_{0}(y,x)$, ${\tau }^{ab}_{i'k'}(y,x;\eta )$, and ${\Phi [h] }^{ab}_{i'k'}(y,x;\eta )$ are evaluated on $(\Sigma ,g(\eta=0))$. 
\end{proposition}

\begin{proof}
From proposition  \ref{usprop2}, taking the limit $\eta \searrow 0^{+}$, we get
\begin{equation}
 g_{i'k'}\,(y,\eta=0)=\left.\int_{\Sigma }{K}^{ab}_{i'k'}(y,x;\eta )\,\left[{g}_{ab}(x,\eta )-2\eta \,\,\mathcal{R}_{ab}(x,\eta )\right]\,d\mu _{g(x,\eta )}\right|_{\forall \,\,\eta>0 }\;,
\end{equation}
which  provides (\ref{grepres}). The asymptotics  follows again from (\ref{uniasi}) under integration term by term and time reversal.
\end{proof}

\noindent Note that explicit expressions for the asymptotic coefficients ${\Phi [h] }^{ab}_{i'k'}(y,x;\beta^{*} )$ can be worked out, at least for the first few terms,  by adapting the relevant formulae in the quoted Gilkey's papers. An interesting application that we will not address here but which seems appropriate to mention at this point is the possibility of (re)-deriving Harnack type estimates, under non-negative curvature assumptions, by directly using the heat kernel of the conjugate Linearized Ricci flow. This application is immediately suggested by the relation (\ref{riccikernel}) and its asymptotics (\ref{ricciasi}) in Proposition \ref{prob1}. 

\noindent The integral representation (\ref{grepres}) of the Ricci flow metric $\beta\mapsto g_{ab}(\beta)$ can be also interpreted as the proof of the following 
\begin{theorem}
The heat kernel flow
\begin{equation}
\eta\longmapsto {K}^{ab}_{i'k'}(y,x;\eta)\;
\end{equation}
is conjugated and thus fully equivalent to the Ricci flow $\beta\longmapsto g_{ab}(\beta)$.
\end{theorem}
\noindent This can be considered as the most important consequence of the conjugacy relation for the linearized Ricci flow. Clearly its utility is somewhat limited by the fact that the flow $\eta\longmapsto {K}^{ab}_{i'k'}(y,x;\eta)$ is constructed on top of the Ricci flow $\beta\mapsto g_{ab}(\beta)$ itself, and thus it does not come as a fully unexpected result. However, it opens to the possibility of a weak formulation of the Ricci flow by exploiting the linear evolution of $\eta\longmapsto {K}^{ab}_{i'k'}(y,x;\eta)$.

\section{Ricci flow conjugated constraint sets}
\noindent To complete the geometrical picture associated with the properties of the conjugate linearized Ricci flow, let us consider, along a Ricci flow of bounded geometry $\beta \mapsto g_{ab}(\beta)$, $0\leq \beta \leq \beta ^{*}$, the heat flow $\beta \longmapsto \varrho (\beta )$, associated with a smooth  function $\varrho (\beta=0)=\rho _{0}\in C^{\infty }(\Sigma ,\mathbb{R})$, \emph{i.e.},
\begin{equation}
\frac{\partial\, \varrho }{\partial \beta}=\Delta \varrho \;,\;\;\;\; \varrho (\beta=0)=\varrho _{0}\;.
\label{scalheat}
\end{equation}
Recall , (see (\ref{adheatop})), that its $L^{2}(M^{4}_{Par})$--conjugate, along the backward Ricci flow $\eta \mapsto g_{ab}(\eta)$, $0\leq \eta \leq \beta ^{*}$, $\eta\doteq \beta^{*}-\beta$, is characterized by the flow $\eta\mapsto \varpi (\eta )$ defined by 
\begin{equation}
\frac{\partial \,\varpi }{\partial \eta}=\Delta \varpi-\mathcal{R}\,\varpi  \;,\;\;\;\; \varpi (\eta=0)=\varpi _{*}\;,
\label{omegadjo}
\end{equation}
where $\varpi _{*}\in C^{\infty }(\Sigma ,\mathbb{R}^{+})$, with $\int_{\Sigma }\,\varpi _{*}\,d\mu_{g(\eta =0)}=1$. 
Since the Riemannian measure is covariantly constant, (\ref{omegadjo}) can be equivalently rewritten as $\frac{\partial  }{\partial \eta}\,d\varpi=\Delta d\varpi$, where $d\varpi (\eta )\doteq \varpi (\eta)\,d\mu_{g(\eta)}$ and $\int_{\Sigma }\,d\varpi(\eta )=1$, $0\leq \eta \leq \beta ^{*}$.
The conjugacy between $\varrho(\beta )$ and $\varpi (\eta)$ is associated with the conservation of the $\varrho(\beta )$--content of $(\Sigma, g_{ab}(\beta ))$ under the flow of probability measures
$\beta \mapsto d\varpi (\eta=\beta ^{*}-\beta)$, \emph{i.e.} 
\begin{equation} 
\frac{d}{d\beta }\,\int_{\Sigma }\,\varrho(\beta )\,d\varpi (\beta)=0\;.
\end{equation}
\noindent The properties of the conjugate heat flow \cite{ecker},\cite{LeiNi},\cite{18} and those of the conjugate linearized Ricci flow established in the previous sections  suggest to shift emphasis from the flows themselves to their dependence from the corresponding initial data.
Thus, along a Ricci flow of bounded geometry $\beta \mapsto (\Sigma ,{g}(\beta ))$, $\beta\in [0,\beta ^{*}]$ let us consider the associated heat flow $(\beta,\varrho _{0}) \mapsto \varrho (\beta )$ and linearized Ricci flow $(\beta,h_{ab}) \mapsto {\widetilde{h}}_{ab}(\beta )$, as functionals of the respective initial data  $\varrho (\beta=0)\doteq \varrho_{0}$, and ${\widetilde h}_{ab}(\beta =0)\doteq h_{ab}$ appearing in the defining PDEs 
(\ref{scalheat}) and (\ref{divfree}), respectively.  In a similar vein let us consider also, along the backward Ricci flow 
$\eta \mapsto (\Sigma ,{g}(\eta ))$, $\eta\in [0,\beta ^{*}]$, $\eta\doteq \beta^{*}-\beta$, the conjugate flows $(\eta,\varpi _{*}) \mapsto \varpi (\eta )$  and $(\eta,H_{*}^{ab}) \mapsto {H}^{ab}(\eta )$, as functionals of the respective initial data  $\varpi (\eta=0)\doteq \varpi_{*}$, and $H^{ab}(\eta =0)\doteq H_{*}^{ab}$ appearing in 
(\ref{omegadjo}) and (\ref{transVol}). 

\noindent For a generic metric $g\in \mathcal{M}et(\Sigma)$, a generic symmetric bilinear form $s_{ik}\in \mathcal{T}_{g}\,\mathcal{M}et(\Sigma )$, and a function $f\in C^{\infty }(\Sigma ,\mathbb{R}^{+})$,  let 
\begin{eqnarray}
\mathcal{T}\,\mathcal{M}et(\Sigma )\times C^{\infty }(\Sigma ,\mathbb{R})&\longrightarrow & \mathbb{R} \\
(g_{ab},s_{ik},f)&\longmapsto & \mathcal{C}(g_{ab},s_{ik},f)=\,0\nonumber\;,
\end{eqnarray}
denote a (surjective) mapping defining a constraint set $\mathcal{C}^{-1}(0)$ in $\mathcal{T}\,\mathcal{M}et(\Sigma )\times C^{\infty }(\Sigma ,\mathbb{R})$, associated with a geometrical condition on the triple $(g_{ab},s_{ik};f)$.

\noindent The following definition is geometrically natural in the light of the properties of the conjugated flows associated with the Ricci flow

\begin{definition} 
Let $\beta \mapsto (\Sigma ,{g}(\beta ))$, $\beta\in [0,\beta ^{*}]$ be a given Ricci flow of bounded geometry, and let $(\beta,\rho _{0}) \longmapsto \varrho (\beta )$ denote the corresponding heat flow associated with the initial condition $\varrho (\beta=0)=\rho _{0}\in C^{\infty }(\Sigma ,\mathbb{R}^{+})$.
If the initial datum ${\widetilde h}_{ab}(\beta =0)\doteq h_{ab}$ for the linearized Ricci flow satisfies the geometrical  constraint
\begin{equation}
\mathcal{C}\left(g_{ab}(\beta =0),h_{ik},\varrho _{0}\right)=\,0\;,
\label{con1}
\end{equation}
and the initial datum $\varpi (\eta=0)\doteq \varpi_{*}$, $H^{ab}(\eta =0)\doteq H_{*}^{ab}$,
for the conjugate heat and the conjugate linearized Ricci flow, can be choosen such that 
\begin{equation}
\mathcal{C}\left(g_{ab}(\eta =0),H_{*}^{ik}, \varpi_{*}\right)=0\;,
\label{con2}
\end{equation}
then the constraints (\ref{con1}) and (\ref{con2}) are said to be conjugated along the given Ricci flow.
\end{definition}
\noindent In order to understand the rationale of such a definition  observe that we cannot expect that a geometrical condition $\mathcal{C}\left(g_{ab}(\beta =0),h_{ik},\varrho _{0}\right)=\,0$ on the initial data will be preserved along  their Ricci flow evolution  $\beta \mapsto (g_{ab}(\beta),{\widetilde h}_{ab}(\beta ), \varrho(\beta ))$. However if, along the associated backward Ricci flow $\eta \mapsto g_{ab}(\eta )$,  we can select initial data $\varpi (\eta=0)\doteq \varpi_{*}$, $H^{ab}(\eta =0)\doteq H_{*}^{ab}$, for the conjugate flows (\ref{omegadjo}) and (\ref{transVol}), such that $\mathcal{C}\left(g_{ab}(\eta =0),H_{*}^{ik},\varpi _{*}\right)=\,0$, then the conjugate flow $\eta \mapsto (g_{ab}(\eta),H^{ab}(\eta ), \varpi(\eta ))$ interpolates between $\mathcal{C}\left(g_{ab}(\beta =0),h_{ik},\varrho _{0}\right)=\,0$ and $\mathcal{C}\left(g_{ab}(\eta =0),H_{*}^{ik}, \varpi_{*}\right)=0$ by averaging the forward flow
$\beta \mapsto (g_{ab}(\beta),{\widetilde h}_{ab}(\beta ), \varrho(\beta ))$ according to the results obtained in section \ref{LRF}, \emph{i.e.}, 
\begin{equation}
\frac{d}{d\beta }\int_{\Sigma }{H}^{ab}(\beta )\,\widetilde{h}_{ab}(\beta )\,d\mu _{g(\beta )}=0\;,
\end{equation}
\begin{equation}
\frac{d}{d\beta }\int_{\Sigma }{H}^{ab}(\beta )\,\mathcal{R}_{ab}(\beta )\,d\mu _{g(\beta )}=0\;,
\end{equation}
\begin{equation}
\frac{d}{d\beta }\int_{\Sigma }{H}^{ab}(\beta )\,{g}_{ab}(\beta )\,d\mu _{g(\beta )}=-2\int_{\Sigma }{H}^{ab}(\beta )\,\mathcal{R}_{ab}(\beta )\,d\mu _{g(\beta )}\;,
\end{equation}
and
\begin{equation} 
\frac{d}{d\beta }\,\int_{\Sigma }\,\varrho(\beta )\,d\varpi (\beta)=0\;.
\end{equation}   

\noindent A typical constraint  $\mathcal{C}\left(g_{ab}(\beta =0),h_{ik},\varrho _{0}\right)=\,0$ one may wish to consider on the triple $\left(g_{ab},h_{ik},\varrho _{0}\right)$ is the of the form
\begin{equation}
\mathcal{R}-|h|^{2}+(tr_{g}(h))^{2}-C\,\varrho_{0} =0\;,
\end{equation}
where we have set $|h|^{2}\doteq h_{ab}h^{ab}$, $tr_{g}(h)\doteq g^{ab}h_{ab}$, and where $C$ is a constant. A constraint of this type occurs in general relativity (the \emph{Hamiltonian constraint}) where it relates the matter density $\varrho_{0}\geq 0$, with the metric $g$, the scalar curvature $\mathcal{R}$ and the second fundamental form $h_{ab}$ of the Riemannian $3$--manifold $(\Sigma ,g)$ carrier of the inital data set for Einstein equations. The above characterization of a Ricci flow conjugated constraint set implies, in this particular setting, that the Hamiltonian constraint is conjugated along a given Ricci flow $\beta \mapsto (\Sigma ,{g}(\beta ))$, $\beta\in [0,\beta ^{*}]$, if we can find triples of initial data $\left(g_{ab}(\beta =0),h_{ik},\varrho _{0}\right)$ and $\left(g_{ab}(\beta =\beta ^{*}),H^{ab}_{*},\varpi _{*}\right)$ such that 
\begin{equation}
\mathcal{R}(\beta =0)-|h|_{\beta =0}^{2}+(tr_{g_{(\beta =0)}}(h))^{2}-C\,\varrho_{0} =0\;,
\end{equation}
and 
\begin{equation}
\mathcal{R(\beta=\beta ^{*})}-|H_{*}|_{\beta =\beta^{*}}^{2}+(tr_{g(\beta=\beta^{*})}(H))^{2}-C\,\varpi_{*} =0\;.
\end{equation}
In such a case, the resulting  conjugate flows $\eta \mapsto (g_{ab}(\eta),H^{ab}(\eta ), \varpi(\eta ))$ interpolates between $\left(g_{ab}(\beta =0),h_{ik},\varrho _{0}\right)$ and $\left(g_{ab}(\eta =0),H_{*}^{ik}, \varpi_{*}\right)$ by formally averaging the Hamiltonian data $\left(g_{ab}(\beta =0),h_{ik},\varrho _{0}\right)$ with the kernels $\left(H^{ab}(\eta), \varpi(\eta)\right)$, \emph{i.e.},
\begin{equation}
\int_{\Sigma }{H}^{ab}(\beta^{*}-\beta)\,\widetilde{h}_{ab}(\beta )\,d\mu _{g(\beta )}=\int_{\Sigma }{H}^{ab}(\beta=0)\,\widetilde{h}_{ab}(\beta=0 )\,d\mu _{g(\beta=0 )}\;,
\end{equation}
\begin{equation}
\int_{\Sigma }{H}^{ab}(\beta^{*}-\beta)\,\mathcal{R}_{ab}(\beta )\,d\mu _{g(\beta )}=\int_{\Sigma }{H}^{ab}(\beta=0)\,\mathcal{R}_{ab}(\beta=0 )\,d\mu _{g(\beta=0 )}\;,
\end{equation}
and
\begin{equation} 
\int_{\Sigma }\,\varrho(\beta )\,d\varpi (\beta-\beta^{*})=\int_{\Sigma }\,\varrho(\beta=0 )\,d\varpi (\beta=0)\;.
\end{equation}   
Thus, we can interpret the existence of a Ricci flow conjugate Hamiltonian constraints as a statement of the possibility of averaging the initial data set $\left(g_{ab}(\beta =0),h_{ik},\varrho _{0}\right)$ over the support of the kernels $\left(H^{ab}(\eta), \varpi(\eta)\right)$. This particular application of the conjugate flows is of potential interest in addressing the possibility of a Ricci flow deformation of initial data sets in General Relativity, and will be discussed in detail elsewhere.
 
\section{Conclusions}
\noindent  The aspects of the conjugated linearized Ricci flow  discussed  here are the most elementary consequences of the conjugacy relation in parabolic spacetime $M^{4}_{Par}$. However, already at this level, they suggests  a number of useful and promising applications to Ricci flow theory. Among these, the study of the stability of singularity formation is perhaps the most interesting. Let us recall that if a solution $\beta\mapsto g_{ab}(\beta)$, $0\leq \beta< T_{0}$, to the Ricci flow develops a singularity at the maximal time $T_{0}$, then such a singularity is said to be a Type--$I$ singularity if $\sup_{\beta\in [0,T_{0})}(T_{0}-\beta)\,\mathcal{K}_{max}(\beta)<+\infty $, whereas it is said to be a Type--$II$ singularity if $\sup_{\beta\in [0,T_{0})}(T_{0}-\beta)\,\mathcal{K}_{max}(\beta)=+\infty $, where $\mathcal{K}_{max}(\beta)$ $\doteq $ $\sup_{x\in \Sigma }\{|\mathcal{R}m\,(x,\beta)|\}$. The analysis of  Type--$II$ singularities is particularly difficult and only recently their existence has been rigorously established for compact manifolds \cite{gu}, (for a nice discussion on Type--$II$ singularities see \cite{6}, \cite{garfisen}, \cite{ding} and \cite{topping}). In particular, since their developments requires a fine tuning between curvature blow--up and neck--pinching, it is not yet clear if they are stable. In known examples, heuristic analysis, and rigorous proofs, Type--$II$ singularities occur when the Ricci flow uses a "critical geometry" for its initial data \cite{garfisen}. Thus, one would expect that a suitable perturbation of such a critical data would remove the degenerate neck--pinching leading to the singularity. However, it is difficult to control what a kind of perturbation would generically remove the criticality. For instance, if $\{x_{i},\beta_{i}\}$ is a sequence of points and of times corresponding to which the curvature, along $\beta\mapsto g_{ab}(\beta)$, attains its maximum, one may think of performing a conformal transformation $\varphi(x_{i},\beta_{i})\, g_{ab}(\beta_{i})$  on the metrics $g_{ab}(\beta_{i})$, and then deform $\varphi(x_{i},\beta_{i})\, g_{ab}(\beta_{i})$  with a corresponding sequence of non--trivial perturbations  $\left\{{\widetilde h}_{ab}(\beta_{i})|\;\;\;\delta\,{\widetilde h}(\beta_{i})=0\right\}$, in such a way that the fine--tuning, between neck--pinching and rounding, leading to the singularity formation, is removed. However, as we have seen, the linearized Ricci flow $\beta\mapsto h_{ab}(\beta)$ does not preserve the non--triviality condition  $\delta\,{\widetilde h}(\beta_{i})=0$, and consequently we do not know a priori which set of perturbations, $(\varphi (\beta=0),h_{ab}(\beta=0))$,  of the critical initial data $g_{ab}(\beta=0)$, will produce the required sequence of deformations $\left\{(\varphi(x_{i},\beta_{i})\, g_{ab}(\beta_{i}),\; {\widetilde h}_{ab}(\beta_{i})  )\right\}$. Thus, the above strategy is difficult, if not impossible, to implement. However, the conjugate linearized Ricci $\eta\mapsto H^{ab}(\eta)$ flow preserves the $\delta\,H=0$ conditions, and one may think to modify the above strategy accordingly: along the sequence $\{x_{i},\beta_{i}\}$ choose conformal factors $\{\phi({i})\}$, and $div$--free $\{H^{ab}({i})\}$ which perturb the sequence of metrics $g_{ab}(\beta_{i})$ by blocking the singularity formation. One can then use the sequence of pairs $\{(\phi({i}),H^{ab}({i}))\}$ as initial data for the conjugate heat flow and for the conjugate linearized Ricci flow. The resulting backward flows $\eta_{i} \longmapsto \{(\phi(\eta_{i}),H^{ab}(\eta_{i}))\}$, with $\eta_{i}\doteq \beta_{i}-\beta$, then generate a sequence of perturbations $\{(\phi(\eta_{i}=\beta_{i}),H^{ab}(\eta_{i}=\beta_{i}))\}$ at $\beta=0$ that can be used to generate, by a limiting procedure, perturbation data $(\phi  ,H^{ab})$ on the initial metric $g_{ab}(\beta=0)$ that will avoid the singularity formation. This is an example where the characterization of Ricci flow conjugated constraint sets appears to be a promising direction for future research.

\vfill\eject

\section*{The heat kernel of a time-dependent generalized Laplacian:\,\,\emph{An appendix by Stefano Romano}}
\label{appendix}
In this appendix we carry out the explicit construction of the the heat kernel of a time-dependent generalized Laplacian, it has been kindly provided by Stefano Romano who adapted to our more general setting the methods used in \cite{guenther2} and in  \cite{chowluni} when treating the scalar time-dependent Laplacian. Although the vector bundle case does not really add anything new from a conceptual point of view, its special importance in the study of the conjugate linearized Ricci flow motivated us to carry out the full computation. Note that here, for technical reasons, we adopt the \emph{analyst sign convention} on Laplacians, \emph{e.g.}, $\Delta :=_{here}-\,g^{ab}\nabla_{ a}\nabla_{ b}$. Also, the result is discussed in the very general setting of vector bundles over a closed manifold carrying a time--dependent metric $g(t)$.

\medskip

\noindent Let $\mathcal{E}\to M^n$ be a vector bundle over a closed manifold $M^n$ and, for $t\in[0, T]$, let $g(t)$ be a time-dependent uniformly bounded family of metrics on $M^n$ and $H_t$ a time-dependent 
family of generalized Laplacians acting on $\Gamma(M^n, \mathcal{E})$. Consider the heat equation
\begin{equation}\left\{
\begin{array}{l}
\displaystyle{\left(\frac{\partial}{\partial t}+H_t\right)s_t=0}\\
\\
s_{t=0}=s_0
\end{array}\right.
\end{equation} 
where $s_t$ is a smooth time-dependent section of $\mathcal{E}$. As usual, $H_t$ determines a unique connection $\nabla_t^\mathcal{E}$ on $\mathcal{E}$ and a unique endomorphism $F_t\in\Gamma(M^n, End(\mathcal{E}))$ such that $H_t=\triangle_t^{\mathcal{E}}+F_t$. We look for a fundamental solution of (3.21), that is a smooth time-dependent section $K_t\in\Gamma(M^n\times M^n, \mathcal{E}\boxtimes\mathcal{E}^*)$ defined for $t>0$ such with the following properties:
\begin{itemize}
\item[(a)] $(\partial_t+H_t)K_t(x, y)=0$, where $H_t$ acts on the $x$ variable, for all $t>0$.
\item[(b)] $\lim_{t\to0}\int_{M^n}K_t(x, y)s(y)d\mu_{g(t)}(y)=s(x)$ for all $s\in\Gamma(M^n, \mathcal{E})$.
\end{itemize}

We refer to condition (b) as the \emph{delta property}.\\
We remark that the notation we have used is imprecise: in fact, since the process we are considering is non-autonomous, i.e. it is not invariant under time-translation, the kernel $K$ should carry explicit dependence on both the inital and final time. By writing $K_t$ we really mean $K_{(t, 0)}$ and we will always use the shorter notation when the initial time is intended to be $t=0$. We will write $K_{(t, \tau)}$ whenever we need to consider a different initial time $\tau\neq0$.\\

Our main result is the following:

\begin{theorem}
\label{theo}
Under the above hypothesis, there exists a unique fundamental solution $K_t(x, y)$ of the heat equation for the time-dependent generalized Laplacian $H_t$.
\end{theorem}

\begin{proof}
We only prove existence, since uniqueness follows easily from properties (a) and (b) by standard methods.\\
We will adopt the technique of first constructing a parametrix for the heat kernel $K_t$ modeled on the euclidean heat kernel and then recovering the full heat kernel form the parametrix. This method only gives a rather generic description of the kernelÕs behavior for small times, but has the advantage of being straightforward. Define 
\begin{equation}
\label{z}
e_t(x, y)\doteq\frac{1}{(4\pi t)^{n/2}}\exp\left(-\frac{d_{g(0)}(x, y)}{4t}\right)
\end{equation}
where $d_{g(0)}$ is the distance function associated to $g(0)$. Choose now a neighborhood $U$ of the diagonal in $M^n\times M^n$ such that $d_{g(0)}(x, y)$ is less than the injectivity radius of $M^n$ for all $(x, y)\in U$. On $U$, define
\begin{equation}
\label{z-}
h^{(K)}_t(x, y)\doteq e_t(x, y)\sum_{\alpha=0}^K\phi_\alpha(x, y; t)t^\alpha
\end{equation}
where both $h^{(K)}_t$ and the $\phi_\alpha$'s are smooth sections of $\mathcal{E}\boxtimes\mathcal{E}^*$ over $U$, and, although $h^{(K)}_t$ is formally defined only for $t>0$, we require the $\phi_\alpha$ to be smooth as $t\to0$. Notice that, differently from the case of a time-independent Laplacian, the inhomogeneity in time of the problem force us to let the $\phi_\alpha$'s depend on time\footnote{This is perhaps clearer if we rewrite (\ref{z-}) with respect to an arbitrary initial time $\tau$ instead of $t=0$:
\begin{equation*}
h^{(K)}_{(t, \tau)}(x, y)=e_t(x, y)\sum_{\alpha=0}^K\phi_\alpha(x, y; \tau)(t-\tau)^\alpha
\end{equation*}
where it is manifest that the kernel $h^{(K)}$ is not invariant under time-translation.}. Our goal is to define the $\phi_\alpha$'s in such a way that
\begin{equation}
\label{a--}
(\partial_t+H_t)h^{(K)}_t=t^Ke_t(\partial_t+H_t)\phi_K
\end{equation}
and that $h^{(K)}_t$ has the delta property.\\
We start by fixing $y$ and choosing coordinates $(x^1, \cdots, x^n)$ near $y$; we also choose normal polar coordinates $(r, \theta^, \cdots, \theta^{n-1})$ centered at $y$ with respect to $g(0)$. Denoting by $\triangle_0$ the scalar Laplacian at $t=0$ and by $J$ the function $\sqrt{\det g(0)}/r^{n-1}$, a standard computation gives
\begin{equation}
\label{a-}
(\partial_t+\triangle_0)e_t=\frac{r}{2t}\frac{\partial \log(J)}{\partial r}e_t
\end{equation}
From now on we drop the subscripts $0$'s and adopt the convention that all quantities that do not exhibit explicit dependence on time refer to $t=0$. We now expand all the relevant quantities in powers of $t$:
\begin{equation}
\label{a}
g_{ij}(t)= g_{ij}+\sum_{\alpha=1}^K(h_\alpha)_{ij}t^\alpha+O(t^{K+1})
\end{equation}
\begin{equation}
\label{b}
\nabla_t^\mathcal{E}= d+\omega(t)=\nabla^\mathcal{E}+\sum_{\alpha=1}^K\omega_\alpha t^\alpha+O(t^{K+1})
\end{equation}
\begin{align}
\label{c}
H_t	&= \triangle_t^\mathcal{E}+F_t=-g_{ij}(t)\Big((\nabla_t^\mathcal{E})_i(\nabla_t^\mathcal{E})_j-\Gamma^k_{ij}(t)(\nabla_t^\mathcal{E})_k\Big)+F_t=\nonumber\\
	&= H+\sum_{\alpha=1}^K\Big(h_\alpha^{ij}\nabla^\mathcal{E}_i\nabla^\mathcal{E}_j+B_\alpha^i\nabla^\mathcal{E}_i+C_\alpha\Big)t^\alpha+ O(t^{K+1})
\end{align}
where the $\omega_\alpha$'s are $End(\mathcal{E})$-valued 1-forms and the $B^i_\alpha, C_\alpha$'s are sections of $End(\mathcal{E})$. Similarly, for the scalar Laplacian we have the expansion
\begin{equation}
\label{d}
\triangle_t=\triangle + \sum_{\alpha=1}^K\left(h^{ij}_\alpha\partial_i\partial_j+b^i_\alpha\partial_i\right) +O(t^{K+1})
\end{equation}
and using (\ref{a-}) we get
\begin{equation*}
(\partial_t+\triangle_t)e_t=\frac{r}{2t}\frac{\partial\log(J)}{\partial r}e_t+\sum_{\alpha=1}^K(h^{ij}_\alpha\partial_i\partial_j+b^i_\alpha\partial_i)t^\alpha e_t+ O(e_t t^{K})=
\end{equation*}
\begin{equation}
\label{e}
=\frac{r}{2t}\Big(\sum_{\alpha=-1}^{K-1}z_\alpha t^\alpha\Big)e_t+ O(e_t t^K)=\frac{r}{2}e_tz_{-1}t^{-1}+\frac{r}{2}e_tz_0+\frac{r}{2}e_t\sum_{\alpha=1}^{K-1}z_\alpha t^\alpha +O(e_tt^K)
\end{equation}
for some smooth functions $z_\alpha$'s. Notice that, since every spatial derivative of $e_t$ brings a factor $t^{-1}$, there is a correction to the lowest order term $\partial_r\log(J)$ coming from the term $th_1^{ij}\partial_i\partial_je_t$.\\
The construction of the parametrix now amounts to expanding everything in $(\partial_t+H_t)h^{(K)}_t$ in powers of $t$, gathering terms of the same order and imposing cancellations up to order $K$. Abbreviating $\Phi^{(K)}_t\doteq \sum_{\alpha=0}^K\phi_\alpha t^\alpha$, we have the formula
\begin{align}
\label{f}
(\partial_t+H_t)h^{(K)}_t= (\partial_t+H_t)e_t\Phi_t^{(K)} 	= &\;((\partial_t+\triangle)e_t)\Phi^{(K)}_t+e_t(\partial_t+H_t)\Phi^{(K)}_t-\nonumber\\
											&-2g^{ij}(t)\partial_ie_t(\nabla_t^\mathcal{E})_j\Phi_t^{(K)}
\end{align}
Using formulas (\ref{a}), (\ref{b}) and (\ref{c}) we find the following expansions:
\begin{align}
\label{g}
(\partial_t+H_t)\Phi^{(K)}_t	&= (\partial_t+H)\phi_0+\phi_1+\nonumber\\
						&+\sum_{\alpha=1}^K\Big((h_\alpha^{ij}\nabla^\mathcal{E}_i\nabla^\mathcal{E}_j+B_\alpha^i\nabla^\mathcal{E}_i+C_\alpha)\phi_0+(\partial_t+H)\phi_\alpha+(\alpha+1)\phi_{\alpha+1}\Big)t^\alpha+\nonumber\\
						&+\sum_{\alpha=2}^K\Big(\substack{\\ \\ \\ \displaystyle{\sum}\\ \gamma, \delta \geq 1\\ \gamma+\delta+\alpha}(h_\gamma^{ij}\nabla^\mathcal{E}_i\nabla^\mathcal{E}_j+B_\gamma^i\nabla^\mathcal{E}_i+C_\gamma)\phi_\delta\Big)t^\alpha +O(t^{K+1})
\end{align}
\begin{align}
\label{h}
-2g^{ij}(t)\partial_ie_t(\nabla_t^\mathcal{E})_j\Phi_t^{(K)}	&= e_tr(\nabla^\mathcal{E}_r\phi_0)t^{-1}+e_tr(\nabla_r^\mathcal{E}\phi_1+(\omega_1)_r\phi_0-h^{ij}_1\partial_ir\nabla^\mathcal{E}_j\phi_0)+\nonumber\\
											&+e_tr\sum_{\alpha=1}^K\bigg[\nabla^\mathcal{E}_r\phi_{\alpha+1}+(\omega_{\alpha+1})_r\phi_0-h^{ij}_{\alpha+1}\partial_ir\nabla_j^\mathcal{E}\phi_0+\nonumber\\
											&+\substack{\\ \\ \\ \displaystyle{\sum}\\ \gamma, \delta\geq1\\\gamma+\delta=\alpha+1}\Big((\omega_\gamma)_r\phi_\delta-h^{ij}_\gamma\partial_ir\nabla^\mathcal{E}_j\phi_\delta-h^{ij}_\gamma\partial_ir(\omega_\delta)_j\phi_0\Big)\bigg]t^\alpha+\nonumber\\
											&+e_tr\sum_{\alpha=2}^K\bigg[\substack{\\ \\ \\ \displaystyle{\sum}\\ \beta, \gamma, \delta\geq1\\ \beta+\gamma+\delta=\alpha+1}-h^{ij}_\beta\partial_ir(\omega_\gamma)_j\phi_\delta\bigg]t^\alpha+ O(e_tt^{K+1})
\end{align}
In the last expression we have exploited the fact that $e_t$ is a radial function in the $g$-normal polar coordinates we have chosen and that $\partial_ie_t=\partial_ir\partial_re_t=(-r/2t)e_t\partial_ir$. We have also written $(\omega_\alpha)_r$ for the radial component of the 1-forms $\omega_\alpha$'s.\\
To finish the computation, we need to substitute the expressions (\ref{e}), (\ref{f}) and (\ref{g}) into (\ref{h}) and impose that all terms of order less than $t^K$ cancel to obtain a chain of differential equations for the $\phi_\alpha$'s. To limit the amount of ugly-looking expressions, we only write down the three lowest order equations:
\begin{equation}
\label{i}
\nabla_r^\mathcal{E}\phi_0+\frac{1}{2}z_{-1}\phi_0=0
\end{equation}
\begin{align}
\label{l}
\nabla_r^\mathcal{E}\phi_1+(\frac{1}{r}+\frac{1}{2}z_{-1})\phi_1 	&=\left(-\frac{1}{2}z_0-\frac{1}{r}(\partial_t+H)-(\omega_1)_r+h^{ij}_1\partial_ir\nabla_j^\mathcal{E}\right)\phi_0=\nonumber\\
													&= F_1(r, \phi_0, z_0; \omega_1, h_1)
\end{align}
\begin{align}
\label{m}
\nabla_r^\mathcal{E}\phi_2+(\frac{2}{r}+\frac{1}{2}z_{-1})\phi_2	&=\bigg(-\frac{1}{2}z_1-\frac{1}{r}(h^{ij}_1\nabla^{\mathcal{E}}_i\nabla^{\mathcal{E}}_j+B^i_1\nabla^{\mathcal{E}}_i+C_1)-\nonumber\\
													&-(\omega_2)_r+h^{ij}_2\partial_ir\nabla_j^{\mathcal{E}}+h_1^{ij}\partial_ir(\omega_1)_j\bigg)\phi_0+\nonumber\\
													&+\bigg(-\frac{1}{2}z_0-(\partial_t+H)-(\omega_1)_r+h_1^{ij}\partial_ir\nabla^{\mathcal{E}}_j\bigg)\phi_1=\nonumber\\
													&=F_2(r, \phi_0, \phi_1, z_0, z_1, \omega_1, \omega_2, h_1, h_2)
\end{align}
and one could continue to arbitrary order to eventually satisfy (\ref{a--}). We still need to verify that $h^{(K)}_t$ has the delta property, but this is easily achieved imposing to equation (\ref{i}) the boundary condition $\phi_0(y, y; 0)= I_{\mathcal{E}_y}$ (in fact, only the 0th order term counts in the small time asymptotic and $e_t(x, y)\to\delta(x, y)$ as $t\to0$).\\
Despite their complicated appearance, equations (\ref{i}), (\ref{l}) and (\ref{m}) are just ODE's to be solved along rays emanating from $y$.\\
Strictly speaking, $h_t^{(K)}$ is not yet the parametrix, since it is only defined in a neighborhood of the diagonal. To extend it to all $M^n\times M^n$, we choose a smooth cutoff function $\eta:[0, \infty)\to[0, 1]$ such that $\eta(x)=1$ if $x\leq\text{inj}(g)/2$ and $\eta(x)=0$ if $x\geq\text{inj}(g)$ and define
\begin{equation}
\label{n}
p_t^{(K)}(x, y)\doteq \eta(d_g(x, y))h_t^{(K)}(x, y)
\end{equation}
for all $(x, y)\in M^n\times M^n$. The last step of the proof consists in constructing the full heat kernel $K_t$ from the parametrix $p_t^{(K)}$. This is achieved by the following
\begin{lemma} Let $p_t^{(K)}$ be a parametrix for the heat equation of a time-dependent Laplacian with $K>n/2$. Then there exists a smooth time-dependent section $\Psi_t\in\Gamma(M^n\times M^n, \mathcal{E}\boxtimes\mathcal{E}^*)$ such that
\begin{equation}
\label{o}
K_t(x, y)\doteq p_t^{(K)}(x, y)+\int_0^td\tau\int_{M^n}p^{(K)}_{(t, \tau)}(x, z)\Psi_\tau(z, y)d\mu_{g(\tau)}(z)
\end{equation}
is the heat kernel.
\end{lemma}

To prove this lemma, we look for $\Psi_t$ as a sum $\sum_{\alpha=1}^\infty(\psi_\alpha)_t$ with the sections $(\psi_\alpha)_t\in\Gamma(M^n\times M^,, \mathcal{E}\boxtimes\mathcal{E}^*)$ defined recursively by
\begin{equation}
\label{p}
(\psi_1)_t(x, y)=(\partial_t+H_t)p^{(K)}_t(x, y)
\end{equation}
and
\begin{equation}
\label{q}
(\psi_{\alpha+1})_t(x, y)=\int_0^td\tau\int_{M^n}\Big[(\partial_t+H_t)p_{(t, \tau)}^{(K)}(x, z)\Big](\psi_{\alpha})_\tau(z, y)d\mu_{g(\tau)}(z)
\end{equation}
Assuming that the series of the $\psi_\alpha$ converges, one checks that the above conditions together with equation (\ref{o}) imply $(\partial_t+H_t)K_t(x, y)=0$. Thus we only need to prove that $\sum_{\alpha=A}^\infty(\psi_\alpha)_t$ converges and that the $K_t$ we have constructed has the delta property.\\
Let $V(t)$ be the volume of $M^n$ at time $t$ and $V\doteq \max_{t\in[0, T]}V(t)$. Fix a fibre metric on $\mathcal{E}$ and define
\begin{equation*}
C=\max_{U\times[0, T]}\Big|(\partial_t+H_t)\phi_K(x, y; t)\Big|
\end{equation*}
with respect to this metric. Then
\begin{equation*}
|(\psi_1)_t|\leq Ct^{K-n/2}
\end{equation*}
and inductively
\begin{equation*}
|(\psi_\alpha)_t|\leq C^\alpha V^{\alpha-1}T^{(\alpha-1)(K-n/2)}\frac{t^{K-n/2+\alpha}}{(K-n/2+\alpha-1)\cdots(K-n/2+1)}
\end{equation*}
and recalling that $K>n/2$ we conclude that $\Psi_t=\sum_{\alpha=1}^\infty(\psi_\alpha)_t$ converges uniformly for all $t$.\\
To see that $K_t$ has the delta property, recall that $p^{(K)}_t$ has it; therefore if we show that the double integral in the right hand side of equation (3.38) vanishes as $t\to0$, we are done. For this condition to be verified, it suffices that $|\int_{M^n}p^{(K)}_{(t, \tau)}(x, z)\Psi_\tau(z, y)d\mu_{g(\tau)}(z)|$ is bounded, or, since $|\Psi|$ is bounded, that $|\int_{M^n}p^{(K)}_{(t, \tau)}(x, z)d\mu_{g(\tau)}(z)|$ is bounded. But this last integral is bounded in the limit $\tau\to t$ because it becomes $|I_{\mathcal{E}_x}|$, as we can easily check using the asymptotics of $p_t^{(K)}$ (recall that we imposed the boundary condition $\phi_0(x, x; 0)=I_{\mathcal{E}_x}$). Moreover, since the metric and the terms $\phi_\alpha$ in the expansion of $p^{(K)}_t$ have uniform bounds over time, we conclude that the integral must be bounded for all $\tau\in[0, T]$. This completes the proofs of the lemma and of Theorem \ref{theo}.
\end{proof}

\bibliographystyle{amsplain}

\end{document}